\documentclass[a4paper,11pt]{article}

\usepackage[latin1]{inputenc}
\usepackage{amsmath,amssymb,amsthm,bbm}

\newtheorem{theorem}{Theorem}[section]
\newtheorem{lemma}[theorem]{Lemma}

\newtheorem{corollary}[theorem]{Corollary}
\newtheorem{definition}[theorem]{Definition}
\newtheorem{example}[theorem]{Example}

\newtheorem{remark}[theorem]{Remark}

\newcommand{\setdef}[2]{\{ #1 \,|\, #2\}}
\newcommand{\bigsetdef}[2]{\bigl\{ #1 \,\bigm|\, #2\bigr\}}

\newcommand{\Def}{\mathrel{\mathop :}=}

\newcommand{\kc}{{\cal C}}

\newcommand{\kh}{{\cal H}}

\newcommand{\kk}{{\cal K}}
\newcommand{\kl}{{\cal L}}

\newcommand{\ko}{{\cal O}}

\newcommand{\IA}{{\mathbb A}}

\newcommand{\IC}{{\mathbb C}}

\newcommand{\IF}{{\mathbb F}}

\newcommand{\IN}{{\mathbb N}}

\newcommand{\IR}{{\mathbb R}}

\newcommand{\IZ}{{\mathbb Z}}

\newcommand{\Waffin}{W^\mathfrak{a}} 
\newcommand{\EWaffin}{\tilde{W}^\mathfrak{a}}
\newcommand{\Haffin}{H^\mathfrak{a}}
\newcommand{\Saffin}{S^\mathfrak{a}}
\newcommand{\EHecke}{\tilde{\kh}^\mathfrak{a}}
\newcommand{\Alko}{\mathcal{A}}
\newcommand{\EAlko}{\tilde{\Alko}}
\newcommand{\Sym}[1]{\mathbf{1}_{#1}}
\newcommand{\Symo}{\mathbf{1}_0}

\newcommand{\Aff}[1]{\langle #1 \rangle_{\text{aff}}}
\newcommand{\Kmg}[1]{\mathcal{#1}}
\newcommand{\Flaff}{\Kmg{G}/\Kmg{B}}

\newcommand{\sfrac}[2]{{\genfrac{}{}{}{3}{#1}{#2}}}
\newcommand{\W}[1]{{W_{#1}(q)}}
\newcommand{\Winv}[1]{{W_{#1}(q^{-1})}}
\newcommand{\pq}{\mathbf{q}}
\newcommand{\w}{\delta}

\numberwithin{equation}{section}

\title{Galleries, Hall-Littlewood polynomials and structure constants of the spherical Hecke algebra}
\author{Christoph Schwer\thanks{This research has been partially supported by the EC TMR network "LieGrits", contract MRTN-CT 2003-505078. \textit{2000 Mathematics Subject Classification. 05E05,33D52}}\\
}

\date{}

\begin{document}

\thispagestyle{empty}

\maketitle

\begin{abstract}
In this paper we provide a combinatorial description of the coefficients appearing in the expansion of Hall-Littlewood polynomials in terms of monomial symmetric functions. We also give a Littlewood-Richardson rule for Hall-Littlewood polynomials. For proving this we use galleries to calculate Satake coefficients and structure constants of spherical Hecke algebras with arbitrary parameters.
\end{abstract}

\section{Introduction}
\label{sec:introduction}

Let $\Phi = (X,\phi,X^\vee,\phi^\vee)$ be a reduced root datum with Weyl group $W$. Let $V = X \otimes \IR$ and $V^* \cong X^\vee \otimes \IR$ such that the pairing $\langle \cdot, \cdot \rangle : V \times V^* \to \IR$ is induced from the perfect pairing between $X$ and $X^\vee$, which will be denoted the same way. Let $Q^\vee$ be the coroot lattice. Choose simple roots $\Delta=\{\alpha_1, \ldots, \alpha_l\}$ and denote by $\phi^+$ the positive roots with respect to $\Delta$. Let $X^\vee_+ = \setdef{x \in X^\vee}{\langle \alpha,x \rangle \geq 0 \text{ for all } \alpha \in \phi^+}$ be the dominant cone. There is a natural action of $W$ on the group algebra $\IZ[X^\vee]$. The algebra of symmetric polynomials $\Lambda$ is the algebra of invariants under this action. It is closely related to the representation theory of complex algebraic groups.

Let $G^\vee$ be a complex reductive linear algebraic group with Borel subgroup $B^\vee$ and maximal torus $T^\vee$ such that the associated root datum is the dual of $\Phi$. Assigning to a finite dimensional representation of $G^\vee$ its $T^\vee$-character yields an isomorphism from the representation ring of $G^\vee$ to~$\Lambda$. For $\lambda \in X^\vee_+$ the Schur polynomial $s_\lambda$ is the character of the irreducible highest weight module $V(\lambda)$ with highest weight $\lambda$. The Schur polynomials $\{s_\lambda\}_{\lambda \in X^\vee_+}$ are a basis of $\Lambda$. The Kostka number $k_{\lambda\mu}$ for $\lambda, \mu \in X^\vee_+$ is the weight multiplicity of $\mu$ in $V(\lambda)$, i.e. the dimension of the $\mu$-weight space $V(\lambda)_\mu$. The Littlewood-Richardson coefficient $c_{\lambda\mu}^\nu$ for $\lambda,\mu,\nu \in X^\vee_+$ is the multiplicity of $V(\nu)$ in $V(\lambda) \otimes V(\mu)$.

Combinatorially the $k_{\lambda\mu}$ and the $c_{\lambda\mu}^\nu$ can be defined as follows: For $\lambda \in X^\vee_+$ define the monomial symmetric function $m_\lambda = \sum_{\mu \in W\lambda} x^\mu$ where for $\mu \in X^\vee$ we denote the corresponding basis element of $\IZ[X^\vee]$ by $x^\mu$. Then $\{m_\lambda\}_{\lambda \in X^\vee_+}$ is a basis of $\Lambda$ and the Kostka numbers are the entries of the transition matrix from the $m_\mu$ to the $s_\lambda$, i.e. we have $s_\lambda = \sum_{\mu \in X^\vee_+} k_{\lambda\mu} m_\mu$. The Littlewood-Richardson coefficients are the structure constants of $\Lambda$ with respect to the Schur polynomials, i.e. $s_\lambda s_\mu = \sum_{\nu \in X^\vee_+} c_{\lambda\mu}^\nu s_\nu$.

One of the main problems of combinatorial representation theory was to give combinatorial formulas for the $k_{\lambda\mu}$ and the $c_{\lambda\mu}^\nu$. One of the solutions is Littelmann's path model in~\cite{littelmann:94}. In~\cite{gaussentlittelmann:03} Gaussent and Littelmann introduced the gallery model and showed that it is equivalent to the path model. They express the $k_{\lambda\mu}$ and the $c_{\lambda\mu}^\nu$ as the number of certain galleries.

In this paper we give analogs of these combinatorial formulas for some $q$-analog of Schur polynomials, the so called Hall-Littlewood polynomials. We describe the transition matrix from monomial symmetric functions to Hall-Littlewood polynomials and we calculate products of Hall-Littlewood polynomials. Specializing $q$ in these formulas we get a new proof for the above mentioned formulas for Schur polynomials in terms of galleries.

Extending the base ring to $\kl^- \Def \IZ[q^{-1}]$ one gets new interesting bases. The Hall-Littlewood polynomials $\{P_{\lambda}(q^{-1})\}$ are a basis for $\Lambda_q \Def \kl^-[X^\vee]^W$ (as $\kl^-$-module). For $\lambda \in X^\vee_+$ they are defined by
\begin{equation*}
P_\lambda(q^{-1}) = \frac{1}{W_\lambda(q^{-1})} \sum_{w \in W} w \Big( x^\lambda \prod_{\alpha \in \phi^+} \frac{1-q^{-1}x^{-\alpha^\vee}}{1-x^{-\alpha^\vee}}\Big).
\end{equation*}
Here $W_\lambda \subset W$ is the stabilizer of $\lambda$ and $W_\lambda(q^{-1}) = \sum_{w \in W_\lambda} q^{-l(w)}$ where $l : W \to \IN$ is the usual length function. The Hall-Littlewood polynomials are $q$-analogs of the Schur polynomials in the sense that $P_\lambda(0) = s_\lambda$. Moreover, we have $P_\lambda(1) = m_\lambda$. For these and other properties of the $P_\lambda(q^{-1})$ see the article~\cite{nelsenram:03} of Nelsen and Ram.

Define Laurent polynomials $L_{\lambda\mu}$ for $\lambda,\mu \in X^\vee_+$ by
\begin{equation*}
P_\lambda(q^{-1}) = \sum_{\mu \in X^\vee_+} q^{-\langle \rho,\lambda + \mu \rangle} L_{\lambda\mu} m_\mu,
\end{equation*}
where $\rho \Def \frac{1}{2} \sum_{\alpha \in \phi^+} \alpha$. By definition we have $q^{-\langle \rho,\lambda + \mu \rangle} L_{\lambda\mu} \in \kl^-$. Moreover, since $P_\lambda(0) = s_\lambda$, the constant term of $q^{-\langle \rho,\lambda + \mu \rangle} L_{\lambda\mu}$ is $k_{\lambda\mu}$. So a combinatorial description of the $L_{\lambda\mu}$ yields a combinatorial description of the $k_{\lambda\mu}$. For non-dominant $\mu \in X^\vee$ we define $L_{\lambda\mu} = q^{\langle \rho, \mu - \mu^+ \rangle} L_{\lambda\mu^+}$, where $\mu^+ \in X^\vee_+$ is the unique dominant element in the $W$-orbit of $\mu$.

The main combinatorial tool for the description of the $L_{\lambda\mu}$ are galleries of generalized alcoves of a fixed type. For details on this and other unexplained notation see section~\ref{sec:galleries}. Following~\cite{gaussentlittelmann:03} we introduce positively folded galleries and associate to each positively folded gallery $\sigma$ a combinatorially defined polynomial $L_\sigma$ in definition~\ref{de:lsigma}. In section~\ref{sec:satakecoefficients} we prove (see also theorem~\ref{th:satakecoefficients} and the paragraph before it)

\begin{theorem}\label{th:basechange}
Let $\lambda \in X^\vee_+$ and $\mu \in X^\vee$. Denote by $t^\lambda$ the type of a minimal gallery from 0 to $\lambda$. Then $L_{\lambda\mu} = q^{-l(w_\lambda)} \sum_{\sigma} q^{l(w_0\iota(\sigma))} L_\sigma$. Here $w_\lambda \in W_\lambda$ is the maximal element and the sum is over all positively folded galleries $\sigma$ of type $t^\lambda$ and weight $\mu$ such that the initial direction $\iota(\sigma)$ is in the set of minimal representatives $W^\lambda$ of $W/W_\lambda$.
\end{theorem}

If $\Phi$ is of type $A$ one gets a combinatorial description of the $L_{\lambda\mu}$ with Young diagrams as a corollary of the results of Haglund, Haiman and Loehr on the monomial expansion of Macdonald polynomials~\cite{haglundhaimanloehr:05}.

We get a description of the $k_{\lambda\mu}$ in terms of galleries by evaluation at $q^{-1} = 0$. We introduce LS-galleries in definition~\ref{def:lsgalleries} (roughly speaking these are the galleries which survive the specialization $q^{-1} = 0$) and get

\begin{corollary}\label{co:kostkanumbers}
For $\lambda, \mu \in X^\vee_+$ the Kostka number $k_{\lambda\mu}$ is the number of LS-galleries of type $t^\lambda$ and weight $\mu$.
\end{corollary}

In definition~\ref{de:lsigma} we also introduce a second monic polynomial $C_\sigma$ for each gallery $\sigma$  which is closely related to $L_\sigma$. We prove that with this statistic one can calculate the structure constants of $\Lambda_q$ with respect to the Hall-Littlewood polynomials. More precisely, define $C_{\lambda\mu}^\nu$ for $\lambda,\mu,\nu \in X^\vee_+$ by
\begin{equation*}
P_\lambda(q^{-1}) P_\mu(q^{-1}) = \sum_{\nu \in X^\vee_+} q^{-\langle \rho, \mu - \lambda + \nu \rangle} C_{\lambda\mu}^\nu P_\nu(q^{-1}).
\end{equation*}
Let $\kc \Def \setdef{x \in V^*}{\langle \alpha,x \rangle \geq 0 \text{ for all } \alpha \in \phi^+}$ be the dominant Weyl chamber. We prove in section~\ref{sec:structureconstants} (see also theorem~\ref{th:structureconstants} and the paragraph before it)

\begin{theorem}\label{th:producthalllittlewood}
Let $\lambda, \mu, \nu \in X^\vee_+$. Then $C_{\lambda\mu}^\nu = q^{-l(w_\mu)} \sum_{\sigma} q^{l(w_0\iota(\sigma))} C_\sigma F_{\mu\nu}^{\varepsilon(\sigma)}$.  Here the sum is over all positively folded galleries $\sigma$ of type $t^\mu$ and weight $\nu$ starting in $\lambda$ such that they are contained in $\kc$ and the final direction $\varepsilon(\sigma)$ is in $W_\nu W^{w_0\mu}$. The correction factor $F_{\mu\nu}^{\varepsilon(\sigma)}$ is contained in $\kl^-$.
\end{theorem}

In~\cite{kapovichmillson:04} a similar formula for $C_{\lambda\mu}^\nu$ with some restrictions on $\lambda,\mu,\nu$ was given by Kapovich and Millson using geodesic triangles in some euclidean building associated to the situation.

For $q^{-1} = 0$ theorem~\ref{th:producthalllittlewood} yields a Littlewood-Richardson rule in terms of galleries.

\begin{corollary}\label{co:littlewoodrichardson}
For $\lambda, \mu, \nu$ in $X^\vee_+$ the Littlewood-Richardson coefficient $c_{\lambda\mu}^\nu$  is the number of LS-galleries $\sigma$ of type $t^\mu$ and weight $\nu - \lambda$ such that the translated gallery $\lambda + \sigma$ is contained in $\kc$.
\end{corollary}

The combinatorial descriptions in the corollaries~\ref{co:kostkanumbers} and~\ref{co:littlewoodrichardson} are more or less the same as the above mentioned descriptions in~\cite{gaussentlittelmann:03}. But our proof is quite different and does not use the combinatorial results stated there. See remark~\ref{re:comparisongaussentlittelmann} for some further details.

For proving theorems~\ref{th:basechange} and~\ref{th:producthalllittlewood} we use the Satake isomorphism to identify $\IZ[q^{\pm\sfrac{1}{2}}][X^\vee]^W$ with the spherical Hecke algebra with equal parameters associated to $\Phi$. Under this isomorphism, Hall-Littlewood polynomials correspond (up to some factor) to the Macdonald basis and the monomial symmetric functions correspond to the monomial basis of the spherical Hecke algebra (see remark~\ref{re:equalparameters}). So the above theorems can be proven in the spherical Hecke algebra.

This is done in the slightly more general setting of spherical Hecke algebras with arbitrary parameters. In section~\ref{sec:galleries} we state our combinatorial formulas regarding Satake coefficients, i.e. the entries of the transition matrix from the monomial basis to the Macdonald basis, and the structure constants of the spherical Hecke algebra with respect to the Macdonald basis. The formulas are then proven in sections~\ref{sec:satakecoefficients} and~\ref{sec:structureconstants}. We define a basis of the affine Hecke algebra indexed by the generalized alcoves introduced in section~\ref{sec:affineweylgroup} and show that right multiplication of this alcove basis by elements of the standard basis can be calculated using galleries.

It is well known that the Satake coefficients form a triangular matrix with respect to the usual partial order on $X^\vee$, i.e. $\mu \leq \lambda$ iff $ \lambda - \mu = \sum_{\alpha \in \Delta} n_\alpha \alpha^\vee$ for nonnegative integers $n_\alpha$. With our combinatorial description of the Satake coefficients we can show that all remaining entries are in fact nonzero.

Now take the spherical Hecke algebra with coefficients in $\IC$ such that all parameters are powers of a fixed prime $p$. Our description of the $L_{\lambda\mu}$ implies that $L_{\lambda\mu} > 0$ for  $\lambda, \mu \in X^\vee_+$ with $\mu \leq \lambda$. This gives a combinatorial proof of a positivity result of Rapoport in~\cite{rapoport:00}.

In section~\ref{sec:commutation} we compute the transition matrix between the alcove basis and the standard basis. This yields a $q$-analog of a commutation formula of Pittie and Ram~(\cite{pittieram:99}) in terms of galleries.

In section~\ref{sec:geometricinterpretation} we give a geometric interpretation of our combinatorics in the case $X^\vee = Q^\vee$. Fix a prime power $\pq$ and let $\IF_\pq$ be the finite field with $\pq$ elements. Let $K$ be the algebraic closure of $\IF_\pq$. Let $G$ be a semisimple, simply connected algebraic group over $K$ with root datum $\Phi$ corresponding to some choice of a Borel $B \subset G$ and a maximal torus $T \subset B$. Assume that all groups are defined and split over $\IF_\pq$.

From the definition of the geometric Satake isomorphism (see for example the survey article~\cite{haineskottwitzprassad:03} of Haines, Kottwitz and Prassad) it is known, that the evaluation at $\pq$ of the coefficients $L_{\lambda\mu}$ gives the number of points of certain intersections in the affine Grassmanian of $G$ over $\IF_\pq$. We show that our combinatorics reflects this interpretation. Using results of Billig and Dyer in~\cite{billigdyer:94}  we show that the galleries occurring in theorem~\ref{th:basechange} together with the associated coefficients parameterize decompositions of these intersections. We also give an interpretation of the alcove basis in this context.

\emph{Acknowledgments.} The author would like to thank P.~Littelmann and A.~Ram for various helpful suggestions and discussions.

\section{Affine Weyl group and alcoves}
\label{sec:affineweylgroup}

In this section we recall some facts on the (extended) affine Weyl group and on alcoves as in~\cite{bourbaki:81}. Furthermore, we introduce generalized alcoves.

The group $Q^\vee$ acts on $V^*$ by translations. The affine Weyl group is defined as the semidirect product $\Waffin = W \ltimes Q^\vee$. It acts on $V^*$ by affine transformations. For $\lambda \in Q^\vee$ denote by $\tau_\lambda \in \Waffin$ the associated translation. The affine Weyl group is generated by its affine reflections. Let $\Haffin$ be the union of all reflection hyperplanes of reflections in $\Waffin$. Then $\Haffin = \bigcup_{\alpha \in \phi^+, m \in \mathbb{Z}} H_{\alpha,m}$, where $H_{\alpha,m} = \setdef{x \in V^*}{\langle \alpha, x \rangle = m}$. Let $H_{\alpha,m}^{\pm} = \setdef{x \in V^*}{\langle \alpha, x \rangle \gtrless m}$ be the associated affine half spaces.

The connected components of $V^* \setminus \Haffin$ are called open alcoves. Their closures are the alcoves in $V^*$. Denote by $\Alko$ the set of all alcoves. The action of $\Waffin$ on $\Alko$ is free and transitive. For $A \in \Alko$ and $\lambda \in Q^\vee$ we have $\tau_\lambda A = \lambda + A = \setdef{\lambda + x}{x \in A}$. 

The fundamental alcove $A_f = \{x \in V^*| 0 \leq \langle \alpha, x \rangle \leq 1 \text{ for all } \alpha \in \phi^+\} \in \Alko$ is a fundamental domain for the $\Waffin$-action on $V^*$. We get a bijection $\Waffin \rightarrow \Alko, w \mapsto A_w \Def wA_f$.

A face $F$ of an alcove $A$ is an intersection $F = A \cap H$ such that $H \subset \Haffin$ is a reflection hyperplane and $\Aff{F} = H$. Here $\Aff{F}$ is the affine subspace spanned by $F$. A wall of $A$ is some hyperplane $H \subset \Haffin$ such that $H \cap A$ is a face of $A$. The group $\Waffin$ is generated by the reflections $\Saffin$ at the walls of $A_f$. One has $\Saffin = S \cup \{s_{01}, \ldots, s_{0c}\}$, where $S = \{s_1, \ldots, s_l\}$ is the set of simple reflections of $W$ and $s_{0k}$ is the affine reflection at $H_{\theta_k,1}$. Here $\setdef{\theta_k}{k = 1, \ldots, c} \subset \phi^+$ is the set of maximal roots with respect to the partial order on $X^\vee$, so $c$ is the number of irreducible components of the Dynkin diagram of $\Phi$. Moreover, $(\Waffin, \Saffin)$ is a Coxeter system.

Let $F$ be a face of $A_f$. The type of $F$ is the reflection at $\Aff{F}$. Extend this definition to all faces by demanding that the $\Waffin$-action preserves types.

Right multiplication of $\Waffin$ induces an action of $\Waffin$ on $\Alko$ from the right. For $A \in \Alko$ and $s \in \Saffin$ the alcove $As$ is the unique alcove not equal to $A$ having a common face of type $s$ with $A$. Let $F_s \subset A$ be the face of type $s$ and $\Aff{F_s} = H_{\alpha,m}$ for some $\alpha \in \phi^+$ and $m \in \IZ$. The hyperplane $H_{\alpha,m}$ is called the separating hyperplane between $A$ and $As$. Call $A$ negative with respect to $s$ if $A$ is contained in $H_{\alpha,m}^-$ and denote this by $A \prec As$. Of course $A$ is called positive with respect to $s$ if $As$ is negative with respect to $s$. We have $A \prec As$ iff $\lambda + A \prec \lambda + As$ for all $\lambda \in Q^\vee$.

\begin{example}\label{ex:alcoveorder}
\begin{itemize}
\item For $A_w$ and $A_{ws}$ in the dominant chamber $\kc$ we have $A_w \prec A_{ws}$ iff $w < ws$, where '$\leq$' is the usual Bruhat order on $\Waffin$.
\item Let $w \in W$ and $s \in S$. Then $A_w \prec A_{ws}$ iff $w > ws$.
\end{itemize}
\end{example}

There is also a natural action of $X^\vee$ on $V^*$ by translations. So we can extend the above definition and get the extended affine Weyl group $\EWaffin \Def W \ltimes X^\vee$. Extending the above notation write $\tau_\mu$ for the translation by $\mu \in X^\vee$. The action of $\EWaffin$ on $\Alko$ is no longer free and type preserving. The stabilizer $\Omega$ of $A_f$ is isomorphic to $X^\vee/Q^\vee$. The isomorphism is given by sending $g \in \Omega$ to the class of $g(0)$. So a set of representatives is given by $X^\vee \cap A_f$. We have $\EWaffin \cong \Omega \ltimes \Waffin$ and every element $v \in \EWaffin$ can be written as $v = wg$ for unique $w \in \Waffin$ and $g \in \Omega$. Although $\EWaffin$ is no longer a Coxeter group, we can extend the definition of the length function by setting $l(v) = l(w)$. So multiplication by elements of $\Omega$ does not change the length. One also can extend the Bruhat order on $\EWaffin$ as follows: Let $v = wg$ and $v' = w'g' \in \EWaffin$ such that $w, w' \in \Waffin$ and $g,g' \in \Omega$. Then $v \leq v'$ iff $g = g'$ and $w \leq w'$ (in the usual Bruhat order on $\Waffin$).

As mentioned above, the action of $\EWaffin$ on $\Alko$ is no longer free. So we introduce generalized alcoves $\EAlko$ in order to work with the extended affine Weyl group as follows: Take an alcove $A \in \Alko$. Then some conjugate of $\Omega$ acts transitively on $A \cap X^\vee$ and this intersection is in natural bijection to $X^\vee / Q^\vee$. Define $\EAlko \Def \setdef{(A,\mu) \in \Alko \times X^\vee}{\mu \in A}$. There is a natural embedding $\Alko \hookrightarrow \EAlko$ sending an alcove $A \in \Alko$ to $(A,\mu)$ where $\mu$ is the unique element in $A \cap Q^\vee$. We identify $\Alko$ with its image in $\EAlko$. We have a natural free $\EWaffin$-action on $\EAlko$ given by the natural action on the two components. In particular, $X^\vee$ acts on $\EAlko$ by translations in both components. For $\lambda \in X^\vee$ and $A \in \EAlko$ we also write $\lambda + A$ for $\tau_\lambda A$.

We get a bijection $\EWaffin \rightarrow \EAlko, w \mapsto A_w \Def w A_f$ extending the bijection $\Waffin \to \Alko$. In the same way as above we also get a right $\EWaffin$-action on $\EAlko$ where $\Omega$ acts only on the second factor. The definitions of face and type of a face carry over to this situation by demanding that $\EWaffin$ acts type preserving.  

Every generalized alcove $A$ is of the form $\mu + A_w$ for unique $\mu \in X^\vee$ and $w \in W$. Then $\mu$ is called the weight of $A$ and $w$ its direction. Denote this by $wt(A) \Def \mu$ and $\w(A) \Def w$.

In various circumstances we will deal with stabilizer subgroups of $W$. We use the following notation for some notions related to them.

\begin{definition}\label{def:stabilizers}
Let $\mu \in X^\vee$ and $W_\mu \subset W$ its stabilizer. The maximal element of $W_\mu$ is denoted by $w_\mu$, the minimal representatives of $W/W_\mu$ by $W^\mu$ and the minimal element in the coset $\tau_\mu W$ by $n^\mu$.
\end{definition}

In particular, $W = W_0$ and $w_0$ is the longest element in $W$. We will frequently use some facts about the length function on $\EWaffin$ summarized in

\begin{lemma}\label{le:lengthfunction}
Let $\lambda \in X^\vee_+$.
\begin{enumerate} 
\item We have $l(\tau_\lambda) = 2 \langle \rho, \lambda \rangle$. In particular, $l$ is additive on $X^\vee_+$.
\item One has $\tau_\lambda w_\lambda = n^\lambda w_0$ and $l(\tau_\lambda) + l(w_\lambda) = l(n^\lambda) + l(w_0)$. Moreover, $n^\lambda \in W \tau_\lambda W$ is minimal.
\end{enumerate}
\end{lemma}

\section{Affine Hecke algebra}
\label{sec:affinehecke}

Details on the affine Hecke algebra of a root datum with unequal parameters can be found in Lusztig's article~\cite{lusztig:89}.

For defining the affine Hecke algebra we first have to fix parameters. Let $d : \Saffin \to \IN$ be invariant under conjugation by elements of $\EWaffin$. Let $\kl = \IZ\big[q^{\pm\sfrac{1}{2}}\big]$ and define $q_s = q^{d(s)}$ for $s \in \Saffin$. For $v \in \Waffin$ we set $q_v = \prod_{j = 1}^k q_{t_{j}}$ where $v = t_{1} \cdot \ldots \cdot t_{k}$ with $t_i \in \Saffin$ is a reduced decomposition of $v$. For arbitrary $v \in \EWaffin$ let $q_v = q_{v'}$ where $v = v'g$ with $v' \in \Waffin$ and $g \in \Omega$.

For a subset $H \subset W$ define $H(q) = \sum_{w \in H} q_w$ and $H(q^{-1}) = \sum_{w \in H} q_w^{-1}$.

The affine Hecke algebra $\EHecke$ associated to the root datum $\Phi$ and the above choice of $d$ is a $\kl$-algebra defined as follows: As a $\kl$-module it is free with basis $\{T_w\}_{w \in \EWaffin}$ and multiplication is given by
\begin{itemize}
\item $T_s^2 = q_s T_{id} + (q_s-1) T_s$ for all $s \in \Saffin$.
\item $T_v T_w = T_{vw}$ for all $v,w \in \EWaffin$ such that $l(vw) = l(v) + l(w)$.
\end{itemize}
On $\EHecke$ there is a natural $\IZ$-algebra involution $\overline{\cdot} : \EHecke \to \EHecke$. It is given by $\overline{T}_w = T_{w^{-1}}^{-1}$ for $w \in \EWaffin$ and $\overline{q^{j}} = q^{-j}$.

For $\lambda \in X^\vee_+$ define $q_\lambda = q^{\sfrac{1}{2}\sum_{j = 1}^k d(t_{j})}$ where $\tau_\lambda = t_1 \cdot \ldots \cdot t_k g$ is a reduced decomposition with $g \in \Omega$. So we have $q_\lambda^2 = q_{\tau_\lambda}$. For arbitrary $\mu \in X^\vee$ define $q_\mu \Def q_\lambda q_{\lambda'}^{-1}$ where $\lambda, \lambda' \in X^\vee_+$ such that $\mu = \lambda - \lambda'$. Then $q_\mu$ is independent of the particular choice of $\lambda, \lambda'$ because of the additivity of the length function on $X^\vee_+$ (see lemma~\ref{le:lengthfunction}).

For each $\mu \in X^\vee$ define an element $X_\mu \in \EHecke$ by $X_\mu \Def q_\mu^{-1} T_{\tau_\lambda} T_{\tau_{\lambda'}}^{-1}$ where as above $\mu = \lambda - \lambda'$ with $\lambda, \lambda' \in X^\vee_+$. By the same reason as above $X_\mu$ does not depend on the choice of $\lambda$ and $\lambda'$. In particular we have $X_\lambda = q_\lambda^{-1} T_{\tau_\lambda}$ for all dominant $\lambda$. Thus one gets an inclusion of $\kl$-algebras 
\begin{align*}
\kl[X^\vee] & \hookrightarrow \EHecke \\
x^\nu & \mapsto X_\nu
\end{align*}
and the image of $\kl[X^\vee]^W$ is the center of $\EHecke$. We identify $\kl[X^\vee]$ with its image.

Define $\Symo = \sum_{w \in W} T_w \in \EHecke$. We have $T_w \Symo = q_w \Symo$ for $w \in W$ and $\Symo^2 = \W{} \Symo$.
The spherical Hecke algebra $\kh^{sph}$ is defined by 
\begin{equation*}
\kh^{sph} = \bigsetdef{h \in \frac{1}{\W{}} \EHecke}{T_w h = h T_w = q_w h \text{ for all } w \in W}.
\end{equation*}
The Macdonald basis of $\kh^{sph}$ is given by $\{M_\lambda\}_{\lambda \in X^\vee_+}$ where
\begin{align*}
M_\lambda & \Def \frac{1}{\W{}} \sum_{w \in W \tau_\lambda W} T_w = \frac{1}{\W{} \W{\lambda}} \Symo T_{n^\lambda} \Symo\\
& = \frac{q_\lambda q_{w_0}^{-1}}{\W{} \Winv{\lambda}} \Symo X_\lambda \Symo.
\end{align*}
The second equality follows from lemma~\ref{le:lengthfunction} which yields $X_\lambda = q_{-\lambda} T_{n^\lambda} T_{w_0} \overline{T}_{w_\lambda}$. Moreover, $\Winv{\lambda} = q_{w_\lambda}^{-1} \W{\lambda}$. One obtains an isomorphism (Satake)
\begin{align*}
\kl[X^\vee]^W & \xrightarrow[]{\cong} \kh^{sph} \\
x & \mapsto \frac{1}{\W{}} x \Symo \\
\intertext{In particular, we have}
m_\lambda & \mapsto Y_\lambda \Def \frac{1}{\W{}} \sum_{\mu \in W\lambda} X_\mu \Symo
\end{align*}

So we have two bases for $\kh^{sph}$: The natural basis given by the Macdonald basis $\{M_\lambda\}_{\lambda \in X^\vee_+}$ and the monomial basis $\{Y_\lambda\}_{\lambda \in X^\vee_+}$ given by the images of the monomial symmetric functions under the Satake isomorphism. We are interested in the transition matrix from the monomial basis to the natural basis. (Re)define the family $\{L_{\lambda\mu}\}_{\lambda, \mu \in X^\vee_+}$ as modified entries of this transition matrix. More precisely, we have \begin{equation*}
M_\lambda = \sum_{\mu \in X^\vee_+} q_\mu^{-1} L_{\lambda\mu} Y_\mu.
\end{equation*}
For arbitrary $\mu \in X^\vee$ and dominant $\lambda \in X^\vee_+$ we set $L_{\lambda\mu} = q_{\mu-\mu^+} L_{\lambda\mu^+}$ where $\mu^+$ is the unique dominant element in the $W$-orbit of $\mu$.

We also calculate the structure constants of the spherical Hecke algebra with respect to the Macdonald basis. For this, (re)define $\{C_{\lambda\mu}^\nu\}_{\lambda,\mu,\nu \in X^\vee_+}$ as modified structure constants by 
\begin{equation*}
M_\lambda M_\mu = \sum_{\nu \in X^\vee_+} q_{\lambda-\nu}^2 C_{\lambda\mu}^\nu M_\nu.
\end{equation*}
In the next chapter we give a description of the $L_{\lambda\mu}$ and the $C_{\lambda\mu}^\nu$ using galleries.

\begin{remark}\label{re:generalparameters}
This is not the most general choice of parameters for which the affine Hecke algebra is defined and where the theorems~\ref{th:satakecoefficients} and~\ref{th:structureconstants} are true. One important example is the following: Replace $\kl$ by the image of the morphism $\kl \to \IC$ evaluating the variable $q$ at some fixed prime power. Hecke algebras of reductive groups over local fields are of this form (compare the end of section~\ref{sec:geometricinterpretation} for the case of equal parameters).
\end{remark}

\begin{remark}\label{re:equalparameters}
Now we want to clarify the relations of this section to symmetric polynomials and their $q$-analogs. In particular, we describe the relation between the coefficients defined above and the ones with the same names in section~\ref{sec:introduction}.

For this regard the case of equal parameters, i.e. $d(s) = 1$ for all $s \in \Saffin$. In this case we have $q_v = q^{l(v)}$ for $v \in \EWaffin$ and $q_\mu = q^{\langle \rho, \mu \rangle}$ for $\mu \in X^\vee$. It is known (see for example~\cite[theorem 2.9]{nelsenram:03}) that the image of $P_\lambda(q^{-1})$ under the Satake isomorphism is $q_{-\lambda} M_\lambda$. So comparing the definitions of the $L_{\lambda\mu}$ and $C_{\lambda\mu}^\nu$ in section~\ref{sec:introduction} with the ones given here shows that the first ones are special cases of the latter ones. So the theorems stated there will follow from theorems~\ref{th:satakecoefficients} and~\ref{th:structureconstants} given in the next section.
\end{remark}

\section{Galleries}
\label{sec:galleries}

In this section we introduce galleries and some polynomials associated to them. We then give a precise meaning to the theorems stated in the introduction in the general setting of the last section. The galleries used here are a slight generalization of the usual galleries in a Coxeter complex since we regard generalized alcoves instead of alcoves.

\begin{definition}\label{de:galleries}
Let $t = (t_1, \ldots, t_k)$ with $t_i \in \Saffin \cup \Omega$. Let $s \in \Saffin$.
\begin{itemize} 
\item 
A gallery $\sigma$ of type $t$ connecting generalized alcoves $A$ and $B$ is a sequence $(A = A_0, \ldots, B = A_k)$ of generalized alcoves such that $A_{i+1} = A_i t_{i+1}$ if $t_{i+1} \in \Omega$ and $A_{i+1} \in \{A_i, A_i t_{i+1}\}$ if $t_{i+1} \in \Saffin$. In the case of $t_{i+1} \in \Saffin$ this means that $A_i$ and $A_{i+1}$ have a common face of type $t_{i+1}$.
\item 
The initial direction $\iota(\sigma)$ is defined to be the direction $\w(A_0)$ of the first generalized alcove. The weight $wt(\sigma)$ of $\sigma$ is $wt(A_k)$, the ending $e(\sigma)$ is $A_k$ and the final direction $\varepsilon(\sigma)$ is $\w(A_k)$.
\item 
The gallery $\sigma$ has a positive $s$-direction at $i$ if $t_{i+1} = s$, $A_{i+1} = A_i s$ and $A_i$ is negative with respect to $s$, i.e. $A_i \prec A_{i+1}$. The separating hyperplane is the wall of $A_i$ corresponding to the face of type $s$.
\item 
The gallery $\sigma$ is $s$-folded at $i$ if $t_{i+1} = s$ and $A_{i+1} = A_i$. The folding hyperplane is the wall of $A_i$ corresponding to the face of type $s$. The folding is positive if $A_i \succ A_i s$.
\end{itemize}
We call $\sigma$ positively folded, if all foldings occurring are positive. A gallery is said to be minimal if it is of minimal length among all galleries connecting the same generalized alcoves.
\end{definition}

For the precise statement on the $L_{\lambda\mu}$ and the $C_{\lambda\mu}^\nu$ we need some statistics on galleries.

\begin{definition}\label{de:lsigma}
Let $\sigma$ be a positively folded gallery of type $t$. For $s \in \Saffin$ define
\begin{itemize}
\item $m_s(\sigma)$ the number of positive $s$-directions.
\item $n_s(\sigma)$ the number of positive $s$-folds.
\item $r_s(\sigma)$ the number of positive $s$-folds such that the folding hyperplane is not a wall of the dominant chamber $\kc$.
\item $p_s(\sigma)$ the number of positive $s$-folds such that the folding hyperplane is a wall of~$\kc$.
\end{itemize}
In particular, $r_s(\sigma) + p_s(\sigma) = n_s(\sigma)$. Now we can define
\begin{itemize}
\item $L_\sigma = \prod_{s \in \Saffin} q_s^{m_s(\sigma)} (q_s-1)^{n_s(\sigma)}$ and
\item $C_\sigma = \prod_{s \in \Saffin} q_s^{m_s(\sigma)+p_s(\sigma)} (q_s-1)^{r_s(\sigma)}$.
\end{itemize}
\end{definition}

For a gallery $\sigma$ such that no folding hyperplane is a wall of $\kc$ one has $L_\sigma = C_\sigma$. In the case of equal parameters (see remark~\ref{re:equalparameters}) we have $\deg L_\sigma = \deg C_\sigma$ for any gallery $\sigma$.

Fix some type $t = (t_1, \ldots, t_k)$. For $A \in \EAlko$ and $\mu \in X^\vee$ let $\Gamma^+_t(A,\mu)$ be the set of all positively folded galleries of type $t$ starting in $A$ with weight $\mu$. Further let $\Gamma^+_t(\mu) = \coprod_{w \in W} \Gamma_t^+(A_w,\mu)$ be the set of all positively folded galleries of weight $\mu$ starting in the origin and let $\Gamma^+_{t}$ be the set of all positively folded galleries starting in the origin. Define
\begin{equation*}
L_t (\mu) \Def \sum_{\sigma \in \Gamma^+_t (\mu)} q_{w_0 \iota(\sigma)} L_\sigma.
\end{equation*}
So there is an additional contribution measuring the distance from $-A_f$ to the initial alcove.

\begin{remark}\label{re:alternativedef}
There is an alternative way of defining $L_t(\mu)$: For any $w \in W$ choose a minimal gallery $\sigma_w$ of type $t_w$ which connects $-A_f$ and $A_w$. Then $\sigma_w$ is a nonfolded gallery of length $l(w_0w) = l(w_0) - l(w)$ and it has only positive directions. The positively folded galleries of type $t_w' = (t_w,t)$ beginning in $-A_f$ correspond to the positively folded galleries of type $t$ starting in $A_{w}$. We get
\begin{equation*}
L_t(\mu) = \sum_{w \in W} \Big( \sum_{\sigma \in \Gamma^+_{t_w'} (-A_f,\mu)} L_\sigma \Big).
\end{equation*}
\end{remark}

Now we can give the formula for the Satake coefficients. Let $\lambda \in X^\vee_+$ and recall the notation introduced in definition~\ref{def:stabilizers}. Let $\sigma^\lambda$ be a minimal gallery connecting $A_f$ and $A_{n^\lambda}$ and denote its type by $t^\lambda$. Using the last definition we get polynomials $L_{t^\lambda}(\mu)$ for all $\mu \in X^\vee$. Up to some factor these are the $L_{\lambda\mu}$. More precisely we prove in section~\ref{sec:satakecoefficients}:

\begin{theorem}\label{th:satakecoefficients}
For $\mu \in X^\vee$ we have
\begin{equation*}
L_{\lambda\mu} = \frac{1}{W_\lambda(q)} L_{t^\lambda} (\mu).
\end{equation*}
Furthermore,
\begin{equation*}
L_{\lambda\mu} = q_{w_\lambda}^{-1} \sum_{\substack{\sigma \in \Gamma^+_{t^\lambda} (\mu)\\ \iota(\sigma) \in W^\lambda}} q_{w_0 \iota(\sigma)} L_\sigma.
\end{equation*}
In particular the $L_{t^\lambda}(\mu)$ do not depend on the choice of the minimal gallery $\sigma^\lambda$ and $L_{t^\lambda} (\mu) = q_{\mu - w\mu} L_{t^\lambda} (w\mu)$ for all $w \in W$.
\end{theorem}

\begin{remark}
One of the surprising implications of the last theorem is the $W$-invariance of the $L_{t^\lambda} (\mu)$ up to some power of $q$. This is surprising because even the cardinality of the sets $\Gamma^+_{t^\lambda} (w\mu)$ depends on $w$.
\end{remark}

\begin{remark}\label{re:minimalgalleries}
Let $w \in \Waffin$. The choice of a minimal gallery $\sigma$ connecting $A_f$ and $A_w$ is equivalent to the choice of a reduced expression for $w$. Let $t = (t_1, \ldots, t_k)$ be the type of $\sigma$. Then we have the reduced expression $w = t_1 \cdot \ldots \cdot t_k$.

Let $v \in \EWaffin$. Then $v$ can be written as $v = wg$ with $w \in \Waffin$ and $g \in \Omega$. A minimal gallery $\sigma$ from $A_f$ to $A_v$ is given by a minimal gallery from $A_f$ to $A_w$ extended by $A_v$. So one can always arrange that at most the last entry of the type of a minimal gallery is in $\Omega$.
\end{remark}

Now it is quite natural to ask when $\Gamma^+_{t^\lambda}(\mu) \neq \emptyset$. Although the definition of galleries is a combinatorial one, it seems hard to give a combinatorial proof for the existence (or non existence) of a gallery of given type and weight. Let $\sigma$ be any gallery of type $t^\lambda$ starting in 0, ending in $A_v$ of weight $\mu$. Since the folding hyperplanes are root hyperplanes we always have $\lambda - \mu \in Q^\vee$. Moreover, $v \leq \iota(\sigma) n^\lambda$ by definition of the Bruhat order on $\EWaffin$. This implies $\mu^+ \leq \lambda$. This also follows from the well known fact that the transition matrix from the monomial basis to the Macdonald basis is triangular with respect to the dominance ordering on~$X^\vee_+$.

The question of the existence of a gallery in $\Gamma^+_{t^\lambda}(\mu)$ does not depend on the choice of parameters $d$. So we can take $d = 1$ as in remark~\ref{re:equalparameters}. Since $P_\lambda(q^{-1})$ and $m_\mu$ are contained in~$\Lambda_q$ we have $q^{-\langle \rho, \lambda + \mu \rangle} L_{\lambda\mu} \in \kl^-$. Moreover, $q^{-l(w_\lambda)} \W{\lambda} = \Winv{\lambda} \in \kl^-$ and thus $q^{-\langle \rho, \lambda + \mu \rangle - l(w_\lambda)} L_{t^\lambda}(\mu) \in \kl^-$. So we get the upper bound
\begin{equation}\label{eq:dimensionestimate}
\deg(L_\sigma) + l(w_0 \iota(\sigma)) \leq \langle \rho, \mu + \lambda \rangle +  l(w_\lambda)
\end{equation}
for all $\sigma \in \Gamma^+_{t^\lambda}(\mu)$. The galleries with maximal degree are of special interest. So define

\begin{definition}\label{def:lsgalleries}
A gallery $\sigma \in \Gamma^+_{t^\lambda}$ is a LS-gallery if we have equality in the above equation, i.e. $\deg(L_\sigma) + l(w_0 \iota(\sigma)) = \langle \rho, wt(\sigma) + \lambda \rangle + l(w_\lambda)$.
\end{definition}

Since $L_\sigma$ is monic we get corollary~\ref{co:kostkanumbers} by evaluating theorem~\ref{th:satakecoefficients} at $q^{-1} = 0$.

The number of LS-galleries in $\Gamma_{t^\lambda}(w\mu)$ is $W$-invariant. This follows from the $W$-invariance (up to a power of $q$) of $L_{t^\lambda}(w\mu)$ in~\ref{th:satakecoefficients}. Now let $\mu \in X^\vee$ such that $\mu^+ \leq \lambda$. We know from representation theory that $k_{\lambda\mu^+} > 0$. So there exists a LS-gallery in $\Gamma^+_{t^\lambda}(\mu^+)$ and thus also in $\Gamma^+_{t^\lambda}(\mu)$ by the $W$-invariance.

Summarizing all this in the following corollary answers the above question on the existence of galleries with a given weight and sharpens the triangularity.

\begin{corollary}\label{co:lsgalleries}
The number of LS-galleries in $\Gamma^+_{t^\lambda} (\mu)$ is $k_{\lambda\mu^+}$. In particular we have $\Gamma^+_{t^\lambda} (\mu) \neq \emptyset$ iff $\mu$ occurs as a weight in $V(\lambda)$, i.e. $\mu^+ \leq \lambda$. Moreover, we have (for arbitrary parameters) $L_{\lambda\mu} \neq 0$ iff $\mu \leq \lambda$.
\end{corollary}

Specializing $q$ at some prime power we get $L_{\lambda\mu} > 0$ for all $\mu \leq \lambda$. This was shown by Rapoport for the case of spherical Hecke algebras of a reductive group over a local field~\cite{rapoport:00}.

\begin{remark}\label{re:comparisongaussentlittelmann}
For regular $\lambda$ the definition of galleries coincides with the one given in~\cite{gaussentlittelmann:03}. Instead of using generalized alcoves they regard galleries of alcoves together with an initial and final weight in $X^\vee$ contained in the first respectively last alcove. This is equivalent to our definition since we can always arrange such that at most the last component of $t^\lambda$ is in $\Omega$ (compare remark~\ref{re:minimalgalleries}). For nonregular $\lambda$ they regard degenerate alcoves. This is more or less the same as our choice of the initial direction. See also remark~\ref{re:minimalrepresentative} for a discussion of this choice.

The proof of corollary~\ref{co:kostkanumbers} in~\cite{gaussentlittelmann:03} is quite different from here. They define root operators on the set of all galleries of type~$t^\lambda$ starting in the origin. Then they show that the subset of LS-galleries is closed under these operators and defines the highest weight crystal with highest weight~$\lambda$.
\end{remark}

We now give the formula for the structure constants replacing $L_\sigma$ with $C_\sigma$. So let $\lambda \in X^\vee_+$ and $t$ be any type. Define $\Gamma^d_{t,\lambda}$ as the set of all positively folded galleries of type $t$ starting in $\lambda$ which are contained in the dominant chamber. Here we allow that folding hyperplanes are contained in the walls of $\mathcal{C}$. For $\nu \in X^\vee_+$ let $\Gamma^d_{t,\lambda}(\nu) \subset \Gamma^d_{t,\lambda}$ be the subset of galleries of weight $\nu$. Define
\begin{equation*}
C_{\lambda t}(\nu) = \sum_{\Gamma^d_{t,\lambda}(\nu)} q_{w_0 \iota(\sigma)} C_\sigma.
\end{equation*}
Now let $\lambda,\mu \in X^\vee_+$ and let $t^\mu$ be the type of a minimal gallery connecting $A_f$ and $A_{n^\mu}$ where $n^\mu \in \tau_\mu W$ is the minimal representative in $\tau_\mu W$. The above definition yields $C_{\lambda t^\mu}(\nu)$ for any $\nu \in X^\vee_+$. Define $F_{\mu\nu}^w \Def q_{w} \sum_{v \in W^{w_0\mu} \cap W_\nu w} q_v^{-1}$ for $\mu, \nu \in X^\vee_+$ and $w \in W$. In section~\ref{sec:structureconstants} we prove:

\begin{theorem}\label{th:structureconstants}
For $\lambda, \mu, \nu \in X^\vee_+$ we have
\begin{equation*}
C_{\lambda\mu}^\nu = \frac{W_\nu(q^{-1})}{W_\mu(q)} C_{\lambda t^\mu}(\nu).
\end{equation*}
Furthermore,
\begin{equation*}
C_{\lambda\mu}^\nu = q_{w_\mu}^{-1} \sum_{\sigma \in \Gamma^d_{t^\mu,\lambda}(\nu)} q_{w_0 \iota(\sigma)} C_\sigma F_{\mu\nu}^{\varepsilon(\sigma)}.
\end{equation*}
In particular, the $C_{\lambda t^\mu}(\nu)$ do not depend on the choice of the minimal gallery.
\end{theorem}

So in contrast to theorem~\ref{th:satakecoefficients} we have a condition on the final direction $\varepsilon(\sigma)$ since $F_{\mu\nu}^w = 0$ iff $w \notin W_\nu W^{w_0\mu}$.

As above we can give an estimate for the degree of the $C_{\lambda t^\nu}(\nu)$ in the case of equal parameters and prove corollary~\ref{co:littlewoodrichardson}. From the last theorem we get $q^{-\langle \rho, \mu - \lambda + \nu \rangle - l(w_\mu)} C_{\lambda t^\mu}(\nu) \in \kl^-$ and thus for any $\sigma \in \Gamma^d_{t^\mu,\lambda}(\nu)$ we have
\begin{equation*}
\deg C_\sigma + l(w_0 \iota(\sigma)) \leq \langle \rho, \mu - \lambda + \nu \rangle + l(w_\mu).
\end{equation*}
Since $\deg L_\sigma = \deg C_\sigma$ and translating a gallery by an element of $X^\vee$ does not change~$L_\sigma$ and the initial direction, corollary~\ref{co:littlewoodrichardson} is proven and we get

\begin{corollary}
For $\lambda, \mu, \nu \in X^\vee_+$ we have $C_{\lambda\mu}^\nu \neq 0$ if $c_{\lambda\mu}^\nu \neq 0 $.
\end{corollary}

Specializing $q$ to a prime power one gets that $C_{\lambda\mu}^\nu > 0$ if $c_{\lambda\mu}^\nu > 0$. For equal parameters this is proven in~\cite{kapovichmillson:04} and also by Haines in~\cite{haines:03} by geometric arguments using the affine Grassmanian of the Langlands dual $G$ of $G^\vee$ to calculate the degree and the leading coefficients of $C_{\lambda\mu}^\nu$.

Another interpretation of the structure constants was given by Parkinson~\cite{parkinson:05} using regular affine buildings.

\section{Satake coefficients}
\label{sec:satakecoefficients}

In this section we introduce the alcove basis of the extended affine Hecke algebra and show that right multiplication of this alcove basis by elements of the standard basis can be calculated using positively folded galleries. From this theorem~\ref{th:satakecoefficients} follows. We also show that one can replace positively folded galleries by negatively folded galleries.

\begin{definition}
Let $A \in \EAlko$. Define $X_A =  q_{-wt(A)} q_{\w(A)} X_{wt(A)} \overline{T}_{\w(A)}$.
\end{definition}

The set $\{X_A\}_{A \in \EAlko}$ is a basis of $\EHecke$. Before we proceed, we need some properties of this basis. First let $\lambda \in X^\vee$ and $A \in \EAlko$. One calculates
\begin{equation}\label{translationalcovebasis}
X_\lambda X_A = q_\lambda X_{\lambda + A}.
\end{equation}
Now assume $A = A_v$ to be dominant such that $\lambda \Def wt(A)$ is regular. Then $v = \tau_\lambda \delta(A)$. Moreover, $\tau_\lambda$ is of maximal length in $\tau_\lambda W$ by lemma~\ref{le:lengthfunction} and $l(v) = l(\tau_\lambda) - l(\w(A))$.  So we get $T_{\tau_\lambda} \overline{T}_{\w(A)} = T_{\tau_\lambda \w(A)} = T_v$ and thus
\begin{equation}\label{dominantalcovebasis}
X_A = q_{-\lambda} q_{\delta(A)} X_\lambda \overline{T}_{\delta(A)} = q_{\tau_\lambda}^{-1} q_{\delta(A)} T_{\tau_\lambda} \overline{T}_{\delta(A)} = q_v^{-1} T_v.
\end{equation}

Multiplying the elements of the alcove basis with $T_s$ from the right can be expressed in terms of the alcove order. It is a $q$-analog of the $\EWaffin$-action on $\EAlko$.

\begin{lemma}\label{le:qmultiplicationalcoves}
Let $A \in \EAlko$. In $\EHecke$ we have
\begin{equation*}
X_A T_s =
\begin{cases}
q_s X_{As} & \text{ if } A \prec As\\
X_{As} + (q_s-1) X_A & \text{ if } A \succ As.
\end{cases}
\end{equation*}
\end{lemma}

\begin{proof}
By~\eqref{translationalcovebasis} the assertion is invariant under translation, i.e. under left multiplication with some $X_\mu$. So it is enough to show the assertion for alcoves $A = A_v$ such that $wt(A)-\alpha^\vee$ is dominant and regular for all $\alpha \in \phi$. By~\eqref{dominantalcovebasis} we have $X_A = q_v^{-1} T_v$ and the multiplication law in $\EHecke$ yields
\begin{equation*}
T_v T_s =
\begin{cases}
T_{vs} & \text{ if } l(v) < l(vs) \\
q_s T_{vs} + (q_s-1) T_{v} & \text{ if } l(v) > l(vs).
\end{cases}
\end{equation*}
But for generalized alcoves in the dominant chamber increasing in the alcove order is equivalent to increasing the length of the corresponding elements of $\EWaffin$ (see example~\ref{ex:alcoveorder}). Moreover, by the choice of $A$ we get $X_{As} = q_{vs}^{-1} T_{vs}$ as elements in $\EHecke$ again by~\eqref{dominantalcovebasis} and the assertion follows.
\end{proof}

Using the same arguments and the fact that multiplying by $T_g$ for $g \in \Omega$ does not change the length we get

\begin{lemma}\label{le:qmultiplicationomega}
For $A \in \EAlko$ we have $X_A T_g = X_{Ag}$ as elements in $\EHecke$.
\end{lemma}

Now we can connect the multiplication in $\EHecke$ to the $L$-polynomials. For generalized alcoves $A$ and $B$ and any type $t$ define $\Gamma^+_t(A,B)$ to be the set of all positively folded galleries of type $t$ connecting $A$ and $B$ and set $L_t(A,B) = \sum_{\sigma \in \Gamma_t^+(A,B)} L_\sigma$.

\begin{lemma}\label{le:recursiongalleries}
Let $t = (t_1, \ldots, t_k), s \in \Saffin$, $t' = (t_1, \ldots, t_k, s)$, and fix generalized alcoves $A$ and $B$. We have
\begin{equation*}
L_{t'}(A,Bs) = 
\begin{cases} 
L_{t}(A,B) & \text{ if } B \succ Bs\\
q_s L_{t}(A,B) + (q_s-1) L_{t}(A,Bs) & \text{ if } B \prec Bs.
\end{cases}
\end{equation*}
\end{lemma} 

\begin{proof}
Let $\sigma' = (A, \ldots, C, Bs) \in \Gamma_{t'}^+(A,Bs)$. Then $C \in \{B, Bs\}$.  We have $C = Bs$ iff $\sigma'$ is $s$-folded at $k+1$. Let $\sigma = (A, \ldots, C)$ and distinguish two cases:\\
$B \succ Bs$: We then have $C = B$ and $\sigma'$ is negative at $k+1$. So $\sigma \in \Gamma_{t}^+(A,B)$ and $L_{\sigma'} = L_\sigma$. Moreover, all galleries in $\Gamma_{t}^+(A,B)$ are obtained this way.\\
$B \prec Bs$: If $C = B$ we have $\sigma \in \Gamma_{t}^+(A,B)$ and $\sigma'$ is positive at $k+1$. So $L_{\sigma'} = q_s L_\sigma$ and one gets all galleries in $\Gamma_{t}^+(A,B)$ this way. If  $C = Bs$ we have $\sigma \in \Gamma_{t}^+(A,Bs)$, $L_{\sigma'} = (q_s-1) L_\sigma$ and one obtains all galleries in $\Gamma_t^+(A,Bs)$ this way.\\
The lemma follows.
\end{proof}

Let $v \in \EWaffin$ and $\sigma$ be a minimal gallery of type $t$ connecting $A_f$ and~$A_v$.

\begin{theorem}\label{th:multiplicationgalleries}
Given $A \in \EAlko$ one has $X_A T_v = \sum_{B \in \EAlko} L_t (A,B) X_B$. 
\end{theorem}

\begin{proof}
Because of lemma~\ref{le:qmultiplicationomega} and since the $L$-polynomials are not affected by elements of $\Omega$ in the type it is enough to show the theorem for $v \in \Waffin$. The proof is done by induction on $l(v)$.\\
Let first $v = s \in \Saffin$: Distinguish two cases.\\
$A \prec As$: In this case $L_{(s)}(A,A) = 0$, $L_{(s)}(A,As) = q_s$ and $L_{(s)} (A,B) = 0$ for all other $B$ and $X_A T_s = q_s X_{As}$.\\
$A \succ As$: In this case $L_{(s)} (A,A) = q_s-1$, $L_{(s)} (A,As) = 1$ and $L_{(s)} (A,B) = 0$ for all other $B$ and $X_A T_s = X_{As} + (q_s-1) X_A$.\\
Now let $v \in \Waffin$, $s \in \Saffin$ such that $l(v) < l(vs)$ and $\sigma' = (A_0, \ldots, A_v, A_{vs})$ is a minimal gallery of type $t'$. Using the last lemma we get
\begin{align*}
X_A T_{vs} & = X_A T_v T_s 
= \Big( \sum_{B \in \Waffin} L_t(A,B) X_B \Big) T_s \\
& = \sum_{B \prec Bs} q_s L_{t} (A,B) X_{Bs} 
+ \sum_{B \succ Bs} L_{t} (A,B) X_{Bs} 
+ \sum_{B \succ Bs} (q_s-1) L_{t} (A,B) X_B\\
& = \sum_{B \prec Bs} \big( q_s L_{t} (A,B) + (q_s-1) L_{t} (A,Bs) \big) X_{Bs} 
+ \sum_{B \succ Bs} L_{t} (A,B) X_{Bs}\\
& = \sum_{B \in \EAlko} L_{t'} (A,Bs) X_{Bs}
= \sum_{B \in \EAlko} L_{t'} (A,B) X_B
\end{align*}
\end{proof}

In particular we get that $L_{t}(A,B)$ does not depend on $\sigma$ and $t$ but only on $v$. Thus the following definition is well defined.

\begin{definition}\label{de:lpolynomials}
For $v \in \EWaffin$ define $L_v (A,B) \Def L_{t}(A,B)$ where $t$ is the type of a minimal gallery from $A_f$ to $A_v$.
\end{definition}

With these results we now can prove proposition~\ref{th:satakecoefficients}.

\begin{lemma}\label{le:macdonaldasmonomial} 
For $\lambda \in X^\vee_+$ we have
\begin{equation*} 
\Symo T_{n^\lambda} \Symo = \sum_{\mu \in X^\vee} q_{-\mu} L_{t^\lambda}(\mu) X_\mu \Symo.
\end{equation*}
\end{lemma}

\begin{proof}
We use the last theorem and the facts that $\overline{\Symo} = q_{w_0}^{-1} \Symo$ and $\overline{T}_w \Symo = q_w^{-1} \Symo$ for all $w \in W$. So one calculates
\begin{align*}
\Symo T_{n^\lambda} \Symo 
& =  q_{w_0} \sum_{w \in W} \overline{T}_w T_{n^\lambda} \Symo 
= q_{w_0} \sum_{w \in W} q^{-1}_w X_{A_w} T_{n^\lambda} \Symo\\
& = q_{w_0} \sum_{w \in W} q_w^{-1} \sum_{\sigma \in \Gamma^+_{t^\lambda}, \iota(\sigma) = w} q_{-wt(\sigma)} q_{\varepsilon(\sigma)} L_\sigma X_{wt(\sigma)} \overline{T}_{\varepsilon(\sigma)} \Symo \\
& = \sum_{w \in W} q_{w_0 w} \sum_{\sigma \in \Gamma^+_{t^\lambda}, \iota(\sigma) = w} q_{-wt(\sigma)} L_\sigma  X_{wt(\sigma)} \Symo\\
& = \sum_{\sigma \in \Gamma^+_{t^\lambda}} q_{w_0 \iota(\sigma)} q_{-wt(\sigma)} L_\sigma X_{wt(\sigma)} \Symo
= \sum_{\mu \in X^\vee} q_{-\mu} L_{t^\lambda} (\mu) X_\mu \Symo
\end{align*}
where the last equality holds by the definition of $L_{t^\lambda}(\mu)$ in section~\ref{sec:galleries}.
\end{proof}

From this we get
\begin{equation*}
M_\lambda = \frac{1}{\W{}\W{\lambda}} \sum_{\mu \in X^\vee} q_{-\mu} L_{t^\lambda}(\mu) X_\mu \Symo.
\end{equation*}
But on the other hand $q_{-\mu} L_{\lambda\mu}$ for dominant $\mu$ is the coefficient of $M_\lambda$ with respect to $Y_\mu$. Moreover, for arbitrary $\nu \in X^\vee$ we defined $L_{\lambda\nu} = q_{\nu - \nu^+} L_{\lambda\nu^+}$. So we get
\begin{align*}
M_\lambda & = \sum_{\mu \in X^\vee_+} q_{-\mu} L_{\lambda\mu} Y_\mu
= \frac{1}{\W{}} \sum_{\mu \in X^\vee_+} \Big( \sum_{\nu \in W\mu} q_{-\nu} L_{\lambda\nu} X_\nu \Symo \Big)\\
& = \frac{1}{\W{}} \sum_{\mu \in X^\vee} q_{-\mu} L_{\lambda\mu} X_\mu \Symo.
\end{align*}
Comparing coefficients of these two expansions we get
\begin{equation*}
L_{\lambda\mu} = \frac{1}{\W{\lambda}} L_{t^\lambda}(\mu)
\end{equation*}
which proves the first statement in~\ref{th:satakecoefficients}. The second statement can be obtained as follows: Every $w \in W$ can be written as $w = w_1 w_2$ for unique $w_1 \in W^\lambda$ and $w_2 \in W_\lambda$ such that $l(w) = l(w_1)+ l(w_2)$ (using the notation introduced in definition~\ref{def:stabilizers}). Define $\Sym{\lambda} = \sum_{w \in W_\lambda} T_w$. Since $\overline{T}_v \overline{T}_w = \overline{T}_{vw}$ for $v, w \in W$ with $l(v)+l(w) = l(vw)$ and $\overline{\Sym{\lambda}} = q_{w_\lambda}^{-1} \Sym{\lambda}$ we get
\begin{equation*}
\Symo = q_{w_0} \sum_{w \in W^\lambda} \overline{T}_{w} \overline{\Sym{\lambda}}
= q_{w_0 w_\lambda} \sum_{w \in W^\lambda} \overline{T}_{w} \Sym{\lambda}.
\end{equation*}
If $v \in W_\lambda$ we have $l(v)+l(n^\lambda) = l(v n^\lambda)$. Moreover, $v n^\lambda = v \tau_\lambda w_\lambda w_0 = \tau_\lambda v w_\lambda w_0 = n^\lambda v'$ with $v' = w_0 w_\lambda v w_\lambda w_0$ by lemma~\ref{le:lengthfunction}. Then $l(v') = l(v)$ and $q_v = q_{v'}$. Thus $T_v T_{n^\lambda} \Symo = T_{n^\lambda} T_{v'} \Symo = q_v T_{n^\lambda} \Symo$ and we get
\begin{equation*}
\Symo T_{n^\lambda} \Symo 
= q_{w_0 w_\lambda} \W{\lambda} \sum_{w \in W^\lambda} \overline{T}_w  \; T_{n^\lambda} \Symo.
\end{equation*}
Now the second statement follows the same way as in the proof of lemma~\ref{le:macdonaldasmonomial} using $q_{w_0 w_\lambda} \sum_{w \in W^\lambda} \overline{T}_w = q_{w_\lambda}^{-1} \sum_{w \in W^\lambda} q_{w_0 w} X_{A_{w}}$.

\begin{remark}\label{re:minimalrepresentative}
In the above considerations there are various other choices for the condition on the initial alcove. Let $v \in W_\lambda$. We have 
\begin{equation*}
\Symo = q_{w_0 w_\lambda v} \sum_{w \in W^\lambda} \overline{T}_{w v} \Sym{\lambda}
\end{equation*} 
since $\overline{T}_v \Sym{\lambda} = q_v^{-1} \Sym{\lambda}$ and the last equation becomes 
\begin{equation*}
\Symo T_{n^\lambda} \Symo = q_{w_0 w_\lambda v} \W{\lambda} \sum_{w \in W^\lambda v} \overline{T}_w  \; T_{n^\lambda} \Symo.
\end{equation*}
Thus one gets
\begin{equation*}
L_{\lambda\mu} = q_{w_\lambda v}^{-1} \sum_{\substack{\sigma \in \Gamma^+_{t^\lambda} (\mu)\\ \iota(\sigma) \in W^\lambda v}} q_{w_0 \iota(\sigma)} L_\sigma.
\end{equation*}
The case considered above was $v = id$. In the case of equal parameters we get for any gallery $\sigma \in \Gamma^+_{t^\lambda}$ such that $\iota(\sigma) \in W^\lambda v$ the upper bound
\begin{equation*}
\deg L_\sigma + l(w_0 \iota(\sigma)) \leq \langle \rho, \lambda + wt(\sigma) \rangle + l(w_\lambda v).
\end{equation*}
One could define LS-galleries to be the ones such that $\iota(\sigma) \in W^\lambda v$ and where there is equality in the last equation. But only with the choice $v = id$ it is enough to impose this equality. The condition on the initial direction follows from this. In particular, for a LS-gallery $\sigma$ we have $\iota(\sigma) \in W^\lambda$.

For the definition of the $L_{t^\lambda}(\mu)$ we started with the minimal representative $n^\lambda$ and we showed that $L_{t^\lambda}(\mu)$ is independent of the initially chosen minimal gallery. One can allow even more freedom in this initial choice. Let $v \in W \tau_\lambda W$ and let $w, w' \in W_\lambda$ such that $v = w n^\lambda w^{'}$ and $l(w) + l(n^\lambda) + l(w') = l(v)$. If instead of $t^\lambda$ we use the type $t$ of a minimal gallery from $A_f$ to $A_v$ we get from the proof of~\ref{le:macdonaldasmonomial} that $L_{t}(\mu) = q_w q_{w'} L_{t^\lambda}(\mu)$ for any $\mu \in X^\vee$. It is clear that the number of LS-galleries in $\Gamma^+_{t}(\mu)$ (with the appropriate changes of the degree condition in the definition) is the same as in $\Gamma^+_{t^\lambda}(\mu)$ since they always encode $s_\lambda$. One also has a canonical bijection between these different sets of LS-galleries. But the total number of galleries in $\Gamma^+_{t}(\mu)$ really depends on the choice of $v$ and this number is minimal if we choose $n^\lambda$. There is another fact that singles out~$n^\lambda$: All the nonfolded galleries are LS-galleries.
\end{remark}

\begin{remark}\label{re:negative}
In definition~\ref{de:galleries} one can replace positive (respectively positively folded) by negative (respectively negatively folded), i.e. one gets $m_s^-(\sigma)$ and $n_s^-(\sigma)$ for each negatively folded gallery $\sigma$. With the obvious changes this yields polynomials $L^-_\sigma$ nonzero only for negatively folded galleries. Going further, one gets $\Gamma^-_t(A,B)$, $L^-_t(A,B)$ and recursions (using the same notations as in~\ref{le:recursiongalleries})
\begin{equation*}
L^-_{t'} (A,Bs) =
\begin{cases} 
L^-_t (A,B) & \text{ if } B \prec Bs\\
q_s L^-_t (A,B) + (q_s-1) L^-_t (A,Bs) & \text{ if } B \succ Bs.
\end{cases}
\end{equation*}
Since $\overline{T}_s = q^{-1}_s (T_s + (1-q_s) T_{id})$ for $s \in \Saffin$ we get from lemma~\ref{le:qmultiplicationalcoves} that
\begin{equation*}
X_A \overline{T}_s =
\begin{cases}
X_{As} + (q^{-1}_s - 1)& \text{ if } A \prec As\\
q^{-1}_s X_{As} & \text{ if } A \succ As
\end{cases}
\end{equation*}
for any $A \in \EAlko$ and $s \in \Saffin$. Under the hypotheses of theorem~\ref{th:multiplicationgalleries} we get
\begin{equation*}
X_A \overline{T}_v = \sum_{B \in \Alko} \overline{L^-_t (A,B)} X_B.
\end{equation*}
If one defines 
\begin{equation*}
L^-_t (\mu) = \sum_{\sigma \in \Gamma^{-}_t (\mu)} q_{\iota(\sigma)} L^-_\sigma
\end{equation*}
we also can express the $L_{\lambda\mu}$ with negatively folded galleries. For this note that left multiplication by $w_0$ on $\EAlko$ induces a type preserving bijection $\phi : \Gamma^+_t \to \Gamma^-_t$ for any type $t$. We have $L^-_{\phi(\sigma)} = L_\sigma$ and $\iota(\phi(\sigma)) = w_0 \iota(\sigma)$. In particular, we get the equality $L_t (\mu) = L^-_t (w_0\mu)$. Combining this with the semi-invariance of the $L_{\lambda\mu}$ with respect to $\mu$ we get
\begin{equation*}
L_{\lambda\mu} = q_{\mu - w_0\mu} L_{\lambda, w_0\mu} = \frac{q^2_\mu}{\W{\lambda}}  L_{t^\lambda}(w_0\mu) = \frac{q^2_\mu}{\W{\lambda}} L^-_{t^\lambda}(\mu)
\end{equation*}
which gives an expression of $L_{\lambda\mu}$ in terms of negatively folded galleries by the definition of $L^-_{t^\lambda} (\mu)$.
\end{remark}

\section{Structure constants}
\label{sec:structureconstants}

In this section we calculate the structure constants of the spherical Hecke algebra with respect to the Macdonald basis and prove theorem~\ref{th:structureconstants} and thus theorem~\ref{th:producthalllittlewood} and its corollary.

\begin{lemma}
Let $A = \mu + A_w$ be a dominant generalized alcove such that $A s$ is no longer dominant. Let $H_{\alpha_i,0}$ be the hyperplane separating $A$ and $A s$. Then we have $X_A T_{s} = T_{s_i} X_A$.
\end{lemma}

\begin{proof}
We have $s_i A = A s$ and $A \succ As$. So $s_i$ and $s$ are conjugate in $\EWaffin$ and thus $q_{s_i} = q_s$. Distinguish two cases:\\
If $s = s_{\theta,1}$ with $\theta \in \Theta$ we have $\langle \alpha_i, \mu \rangle = 1$ and thus $s_i(\mu) = \mu - \alpha_i^\vee$ and $s_i A = s_i(\mu) + A_{s_i w}$. But on the other hand we have $A s = (\mu + w\theta_k^\vee) + A_{w s_\theta}$ and so $w\theta^\vee = -\alpha_i^\vee$. In particular, $s_i w < w$. From~\cite[lemma 2.7(d) and proposition 3.6]{lusztig:89} we know that $q_{\alpha_i^\vee} = q_s$ in this case and
\begin{equation*}
T_{s_i} X_\mu = X_{\mu -\alpha_i^\vee} T_{s_i} + (q_{s_i} - 1) X_\mu.
\end{equation*}
Together with $s_i w < w$ this yields
\begin{equation*}
T_{s_i} X_\mu \overline{T}_w = X_{\mu-\alpha_i^\vee} \overline{T}_{s_i w} + (q_{s_i} - 1) X_\mu \overline{T}_w
\end{equation*}
and thus 
\begin{equation*}
T_{s_i} X_A = X_{A_s} + (q_{s_i} -1) X_A = X_A T_s
\end{equation*}
where the last equality follows from $A \succ As$.\\
If $s = s_j \in S$ we have $s_i(\mu) = \mu$ and $w^{-1}(\alpha_i) = \alpha_j$. So here $s_i w > w$. Using $T_{s_i} X_\mu = X_\mu T_{s_i}$ one obtains the desired equality as above.
\end{proof}

We keep the notation of the last lemma and get $\Symo X_A T_{s}  = \Symo T_{s_i} X_A  = q_s \Symo X_A$ (recall that $q_{s_i} = q_s$). For a generalized alcove $A$ and a type $t$ define $\Gamma^+_{t,A}$ to be the set of all positively folded galleries of type $t$ with initial alcove $A$.

Let $t = (t_1, \ldots, t_k)$ be a type and define $T_t = T_{t_1} \cdot \ldots \cdot T_{t_k}$. From theorem~\ref{th:multiplicationgalleries} we get
\begin{equation*}
X_A T_t = \sum_{\sigma \in \Gamma^+_{t,A}} L_\sigma X_{e(\sigma)}
\end{equation*}
where $e(\sigma)$ is the ending of $\sigma$ as introduced in definition~\ref{de:galleries}. This yields
\begin{equation*}
\Symo X_A T_t = \sum_{\sigma \in \Gamma^+_{t,A}} L_\sigma \Symo X_{e(\sigma)}.
\end{equation*}
Setting $t' = (s, t)$ we obtain by the same arguments
\begin{equation*}
\Symo X_A T_{s} T_t = \sum_{\sigma \in \Gamma^+_{t',A}} L_\sigma \Symo X_{e(\sigma)}.
\end{equation*}
Since $\Symo X_A T_s T_t= q_s \Symo X_A T_t$ we get the following

\begin{lemma}\label{le:latwalls1}
Let $t$ be any type and let $A$ be a dominant generalized alcove such that $As$ is no longer dominant. Setting $t' = (s,t)$ we have
\begin{equation*}
q_s \sum_{\sigma \in \Gamma^+_{t,A}} L_\sigma \Symo X_{e(\sigma)}
= \sum_{\sigma \in \Gamma^+_{t',A}} L_\sigma \Symo X_{e(\sigma)}.
\end{equation*}
\end{lemma}

Now let $\lambda \in X^\vee_+$. Then the generalized alcove $A \Def \lambda + A_w$ is dominant iff $w^{-1} \in W^\lambda$. Let $w^{-1} \in W^\lambda$ and $v \in W_\lambda$. Since $\overline{T}_v X_\lambda = X_\lambda \overline{T}_v$ we get
\begin{equation*}
\overline{T}_v X_{A} = q^{-1}_{v} X_{\lambda+A_{vw}} = q^{-1}_{v} X_{vA}.
\end{equation*}
Since $v \in W_\lambda$ we get the equality (using the notation introduced before the last lemma)
\begin{equation*}
\Symo X_{vA} T_t = q_v \Symo \overline{T}_{v} X_A T_t = \Symo X_{A} T_t.
\end{equation*}
For later use observe that $v A = \lambda + A_{vw}$ and thus $v A$ is no longer dominant. We get

\begin{lemma}\label{le:latwalls2}
Let $\lambda \in X^\vee_+$, $w^{-1} \in W^\lambda$ and $v \in W_\lambda$. Let $A = \lambda + A_w$. For any type $t$ we have
\begin{equation*}
\sum_{\sigma \in \Gamma^+_{t,A}} L_\sigma \Symo X_{e(\sigma)}
= \sum_{\sigma \in \Gamma^+_{t,vA}} L_\sigma \Symo X_{e(\sigma)}.
\end{equation*}
\end{lemma}

Now let $\lambda, \mu \in X^\vee_+$. Let $w_\mu \in W_\mu$ and $n^\mu \in \tau_\mu W$ as in definition~\ref{def:stabilizers}. Let $t^\mu$ denote the type of a minimal gallery from $A_f$ to $A_{n^\mu}$. As in the proof of lemma~\ref{le:macdonaldasmonomial} we get
\begin{align}
\Symo X_\lambda \; \Symo T_{n^\mu} 
& = \Symo X_\lambda \sum_{\sigma \in \Gamma^+_{t^\mu}} q_{w_0 \iota(\sigma)} L_\sigma X_{e(\sigma)}\\
\label{eq:multiplication}
& = q_{\lambda} \sum_{\sigma \in \Gamma^+_{t^\mu,\lambda}} q_{w_0 \iota(\sigma)} L_{\sigma} \Symo X_{e(\sigma)}.
\end{align}
Here $\Gamma^+_{t^\mu,\lambda}$ is the set of all galleries of type $t^\mu$ starting in $\lambda$ and the last equality holds since translating a gallery $\sigma$ by $\lambda$ does not change $L_\sigma$. So we have an expansion for the product in terms of $X_{A}$ for $A \in \EAlko$. But we need the expansion in terms of $X_A$ for dominant $A$ to compute the structure constants.

\begin{theorem}\label{th:reducingtodominant}
For $\lambda, \mu \in X^\vee_+$ we have
\begin{equation*} 
\Symo X_\lambda \Symo T_{n^\mu}
= q_{\lambda} \Winv{\lambda} \sum_{\sigma \in \Gamma^d_{t^\mu,\lambda}} q_{w_0\iota(\sigma)} C_\sigma \Symo X_{e(\sigma)}.
\end{equation*}
\end{theorem}

\begin{proof}
For the proof of this theorem we use lemmas~\ref{le:latwalls1} and~\ref{le:latwalls2} to show that the contribution of the galleries with non-dominant weights in the formula~\eqref{eq:multiplication} is exactly the contribution of the $p_s$.\\
First assume $\lambda$ is regular. Then the first generalized alcove of every gallery starting in $\lambda$ is dominant. Let $\eta \in \Gamma^+_{t^\mu,\lambda}$ be a gallery leaving the dominant chamber. Let $\gamma$ be the maximal initial subgallery of $\eta$ contained in $\kc$ and let $A$ be $e(\gamma)$. Then $\eta$ is not folded after $A$ and the next generalized alcove in $\eta$ is of the form $As$ for some $s \in \Saffin$. Denote by $\Gamma^+_{\gamma} \subset \Gamma^+_{t^\mu,\lambda}$ the set of galleries starting with $\gamma$. By lemma~\ref{le:latwalls1} we have that
\begin{equation*}
\frac{q_s}{q_s-1} \sum_{\sigma \in \Gamma^+_{\gamma}, \sigma \text{ folded at } A} L_\sigma \Symo X_{e(\sigma)}
= \sum_{\sigma \in \Gamma^+_{\gamma}}  L_\sigma \Symo X_{e(\sigma)}.
\end{equation*}
So the contribution of all galleries starting with $\gamma$ is the same as the contribution of the galleries starting with $\gamma$ and staying in $\kc$ at $A$, if the contribution of the folding at $A$ is $q_s$ instead of $q_s-1$. Iteration of this procedure eventually yields
\begin{equation*}
\sum_{\sigma \in \Gamma^+_{\gamma}} L_\sigma \Symo X_{e(\sigma)}
= \sum_{\sigma \in \Gamma^+_{\gamma}, \sigma \subset \kc} C_\sigma \Symo X_{e(\sigma)}
\end{equation*}
which proves the theorem for regular $\lambda$.\\
If $\lambda$ is non-regular we have to apply lemma~\ref{le:latwalls2} to obtain the theorem because in this case the first alcove of a gallery starting in $\lambda$ can be non-dominant. In this case its contribution has a part coming from the initial direction, which we did not need to consider in the regular case. But lemma~\ref{le:latwalls2} tells us that the contribution arising from these alcoves is the same as the contribution from the dominant ones. More precisely we have for $w^{-1} \in W^\lambda$ and $v \in W_\lambda$
\begin{equation*}
\sum_{\sigma \in \Gamma^+_{t^\mu,\lambda}, \iota(\sigma) = w} q_{w_0w} L_\sigma \Symo X_{e(\sigma)}
= q_v \sum_{\sigma \in \Gamma^+_{t^\mu,\lambda}, \iota(\sigma) = vw} q_{w_0vw} L_\sigma \Symo X_{e(\sigma)}
\end{equation*}
and thus
\begin{equation*}
\Winv{\lambda} \sum_{\sigma \in \Gamma^+_{t^\mu,\lambda}, \iota(\sigma) = w} q_{w_0w} L_\sigma \Symo X_{e(\sigma)}
= \sum_{\sigma \in \Gamma^+_{t^\mu,\lambda}, \iota(\sigma) \in W_\lambda w} q_{w_0\iota(\sigma)} \Symo X_{e(\sigma)}.
\end{equation*}
Since the sum over all $w^{-1} \in W^\lambda$ of the left hand side of the last equation is exactly the contribution of the galleries starting in $\kc$, the theorem follows.
\end{proof}

\begin{remark}
The proofs for multiplying Schur polynomials using paths are of a similar type as above (see for example~\cite[section~6]{littelmann:94}). First one gets a formula involving also Schur polynomials associated to paths leaving the dominant chamber. Then one shows that the contributions of the leaving paths cancel each other. This is done by combinatorial arguments, i.e. one can see which paths cancel each other. In contrast to this we do not have any concrete information about this cancellation process.
\end{remark}

Now we can prove the first part of theorem~\ref{th:structureconstants} respectively theorem~\ref{th:producthalllittlewood}. We multiply the equation of the last theorem from the right by $\Symo$ and get by the definition of the Macdonald basis
\begin{align*}
M_\lambda M_\mu
&= \frac{q_{\lambda} q_{w_0}^{-1}}{W(q)\Winv{\lambda}} \frac{1}{\W{}\W{\mu}} \Symo X_\lambda \Symo \, \Symo T_{n^\mu} \Symo\\
&= \frac{q_{\lambda}^2 q_{w_0}^{-1}}{\W{} \W{\mu}} \sum_{\sigma \in \Gamma^d_{t^\mu,\lambda}} q_{w_0\iota(\sigma)} C_\sigma \Symo X_{e(\sigma)} \Symo\\
&= \frac{q_{\lambda}^2 q_{w_0}^{-1}}{\W{} \W{\mu}} \sum_{\sigma \in \Gamma^d_{t^\mu,\lambda}} q_{-wt(\sigma)} q_{w_0\iota(\sigma)} C_\sigma \Symo X_{wt(\sigma)} \Symo\\
&= \frac{q_{\lambda}^2}{\W{\mu}} \sum_{\sigma \in \Gamma^d_{t^\mu,\lambda}} q_{-wt(\sigma)}^2 q_{w_0\iota(\sigma)} C_\sigma \Winv{wt(\sigma)} M_{wt(\sigma)}\\
& = \frac{q_{\lambda}^2}{\W{\mu}} \sum_{\nu \in X^\vee_+} q_{-\nu}^2 \Winv{\nu} C_{t^\mu,\lambda}(\nu) M_\nu.
\end{align*}

To prove the second part of theorem~\ref{th:structureconstants} and thus theorem~\ref{th:producthalllittlewood} we need one more step. It is not possible to impose conditions on the initial direction as in theorem~\ref{th:satakecoefficients}. Instead we impose conditions on the final direction to get rid of the fraction $\frac{1}{\W{\mu}}$. For doing this we need some preparation. The situation is more difficult than the case of Satake coefficients since now two stabilizers instead of one are involved. So we first need some information on the interplay between them.

We use the notation for stabilizer subgroups introduced in definition~\ref{def:stabilizers}. Moreover, for any $\nu \in X^\vee$ let $\Sym{\nu} = \sum_{w \in W_\nu} T_w$ be the corresponding symmetrizer. Note that $W_{w_0\mu} = w_0 W_{\mu} w_0$ and thus $q_{w_\mu} = q_{w_{w_0\mu}}$ and $\W{\mu} = \W{w_0\mu}$.

Let $Y = \sum_{w \in W} R_w \overline{T}_w \in \EHecke$ with $R_w \in \kl$. Assume $Y \in \EHecke \Sym{w_0\mu}$. Then $R_{w} = R_{wv}$ for any $w \in W$ and $v \in W_{w_0\mu}$ and thus
\begin{equation}\label{eq:endingdirection}
Y = q_{w_\mu}^{-1} \sum_{w \in W^{w_0\mu}} R_w \overline{T}_{w} \Sym{w_0\mu}
\end{equation}
since for $w \in W^{w_0\mu}$ we have $\overline{T}_w \overline{\Sym{w_0\mu}} = \sum_{v \in W_{w_0\mu}} \overline{T}_{w v}$ and $\overline{\Sym{w_0 \mu}} = q_{w_\mu}^{-1} \Sym{w_0\mu}$.

Let $\nu \in X^\vee_+$ and take $Y$ of a special form, namely $Y = \sum_{w^{-1} \in W^\nu} R_w \Sym{\nu} \overline{T}_w$. For $w \in W$ denote by $w^\nu$ the minimal element of the coset $W_\nu w$. In particular $(w^\nu)^{-1} \in W^\nu$. Expanding $Y$ in terms of the $\overline{T}_w$ yields
\begin{equation*}
Y = q_{w_\nu} \sum_{w \in W} R_{w^\nu} \overline{T}_w.
\end{equation*}
So if in addition $Y \in \EHecke \Sym{w_0\mu}$ we get $Y = q_{w_\nu} q_{w_\mu}^{-1} \sum_{w \in W^{w_0\mu}} R_{w^\nu} \overline{T}_{w} \Sym{w_0\mu}$ by the considerations above.

We calculate $Y \Symo$ and get $Y \Symo = q_{w_\nu} q_{w_\mu}^{-1} \W{\mu} \sum_{w \in W^{w_0\mu}} q_w^{-1} R_{w^\nu} \Symo$. Thus
\begin{align}\label{eq:doublecosetmultiplication}
Y \Symo = q_{w_\nu} \Winv{\mu} \sum_{w^{-1} \in W^\nu} q_w^{-1} F_{\mu\nu}^{w} R_{w} \Symo
\end{align}
where $F_{\mu\nu}^{w} \Def q_w \sum_{v \in W^{w_0\mu} \cap W_\nu w} q_v^{-1}$. Observe that $W^{w_0\mu} \cap W_\nu w \neq \emptyset$ iff $w \in W_\nu W^{w_0\mu}$. In particular, we get for regular $\nu$ that $F_{\mu\nu}^{w} = 1$ if $w \in W^{w_0\mu}$ and 0 else.

Now we relate this to our problem. We have
\begin{equation*}
\W{\mu} \Symo T_{n^\mu} = \Symo \Sym{\mu} T_{\tau_\mu} T_{w_\mu}  \overline{T}_{w_0} 
= \Symo T_{\tau_\mu} T_{w_\mu} \Sym{\mu} \overline{T}_{w_0} = \Symo T_{n^\mu} T_{w_0} \Sym{\mu} \overline{T}_{w_0}.
\end{equation*}
But $T_{w_0} T_w \overline{T}_{w_0} = T_{w_0 w w_0}$ for all $w \in W$ and thus $T_{w_0} \Sym{\mu} \overline{T}_{w_0} = \Sym{w_0\mu}$. So $\Symo T_{n^\mu} \in \EHecke \Sym{w_0\mu}$.

We have $\W{\nu} \Symo X_\nu = \Symo X_\nu \Sym{\nu}$ since $T_w X_\nu = X_\nu T_w$ for any $w \in W_\nu$. So
the contribution of $\Symo X_\nu$ in theorem~\ref{th:reducingtodominant} is given by
\begin{equation*}
\sum_{\sigma \in \Gamma^d_{t^\mu,\lambda}(\nu)} q_{w_0\iota(\sigma)} C_\sigma \Symo X_{e(\sigma)}
= \frac{q_{-\nu}}{\W{\nu}} \Symo X_\nu \sum_{\sigma \in \Gamma^d_{t^\mu,\lambda}(\nu)} q_{w_0\iota(\sigma)} q_{\varepsilon(\sigma)} C_\sigma \Sym{\nu} \overline{T}_{\varepsilon(\sigma)}.
\end{equation*}
As already observed before, $\nu + A_v \subset \kc$ with $v \in W$ iff $v^{-1} \in W^\nu$. So the final directions of the galleries $\sigma$ occurring in the last equation satisfy $(\varepsilon(\sigma))^{-1} \in W^\nu$. If we define $Y \Def \sum_{\sigma \in \Gamma^d_{t^\mu,\lambda}(\nu)} q_{w_0\iota(\sigma)} q_{\varepsilon(\sigma)} C_\sigma \Sym{\nu} \overline{T}_{\varepsilon(\sigma)}$ then $Y$ is of the kind considered above and $Y \in \EHecke \Sym{w_0\mu}$. So we can apply~\eqref{eq:doublecosetmultiplication} and get
\begin{equation*}
\sum_{\sigma \in \Gamma^d_{t^\mu,\lambda}(\nu)} q_{w_0\iota(\sigma)} q_{\varepsilon(\sigma)} C_\sigma \Sym{\nu} \overline{T}_{\varepsilon(\sigma)} \Symo
= q_{w_\nu} \Winv{\mu} \sum_{\sigma \in \Gamma^d_{t^\mu,\lambda}(\nu)} q_{w_0 \iota(\sigma)} C_\sigma F_{\mu \nu}^{\varepsilon(\sigma)} \Symo.
\end{equation*}
Bringing all this together we can multiply the assertion of theorem~\ref{th:reducingtodominant} from the right by $\Symo$ and get
\begin{equation*}
\Symo X_\lambda \Symo T_{n^\mu} \Symo 
= q_\lambda \Winv{\lambda} \Winv{\mu} \sum_{\nu \in X^\vee_+} \frac{q_{-\nu}}{\Winv{\nu}} \sum_{\sigma \in \Gamma^d_{t^\mu,\lambda}(\nu)} q_{w_0 \iota(\sigma)} C_\sigma F_{\mu \nu}^{\varepsilon(\sigma)} \Symo X_\nu \Symo.
\end{equation*}

Now we can calculate the coefficient of $M_\nu$ in the product $M_\lambda M_\mu$ as above. It is equal to
\begin{equation*}
q_{\lambda-\nu}^2 q_{w_\mu}^{-1} \sum_{\sigma \in \Gamma^d_{t^\mu,\lambda}(\nu)} q_{w_0 \iota(\sigma)} C_\sigma F_{\mu \nu}^{\varepsilon(\sigma)}
\end{equation*}
which proves the second part of theorem~\ref{th:structureconstants} and thus~\ref{th:producthalllittlewood}.

\begin{remark}
Consider the case of equal parameters and let $w^{-1} \in W^\nu$. Then we have $F_{\mu\nu}^w = q^{l(w)} \sum_{v \in W^{w_0\mu} \cap W_\nu w} q^{-l(v)}$. By definition of $W^\nu$ we have $l(v) \geq l(w)$ for all $v \in W_\nu w$ and thus $F_{\mu\nu}^w \in \kl^-$. Moreover, the constant term of $F_{\mu\nu}^w$ is 1 iff $w \in W^{w_0\mu}$.
\end{remark}

\begin{remark}
One can proceed the same way to obtain a formula for the Satake coefficients as in the second part of theorem~\ref{th:satakecoefficients} with a condition on the final direction. For stating the results we consider again the situation of section~\ref{sec:satakecoefficients}. So $\lambda \in X^\vee_+$ and $t^\lambda$ is the type of a minimal gallery from $A_f$ to $A_{n^\lambda}$. Applying the above considerations (for $\lambda$ instead of $\mu$) yields $\Symo T_{n^\lambda} \in \EHecke \Sym{w_0\lambda}$. A formula for $\Symo T_{n^\lambda}$ is given by (see the proof of lemma~\ref{le:macdonaldasmonomial}) $\sum_{\sigma \in \Gamma_{t^\lambda}^+} q_{-wt(\sigma)} q_{\varepsilon(\sigma)} q_{w_0\iota(\sigma)} L_\sigma X_{wt(\sigma)} \overline{T}_{\varepsilon(\sigma)}$. So we get by~\eqref{eq:endingdirection}
\begin{equation*}
\Symo T_{n^\lambda} = q_{w_\lambda}^{-1} \sum_{\substack{\sigma \in \Gamma_{t^\lambda}^+\\ \varepsilon(\sigma) \in W^{w_0\lambda}}} q_{-wt(\sigma)} q_{\varepsilon(\sigma)} q_{w_0\iota(\sigma)} L_\sigma X_{wt(\sigma)} \overline{T}_{\varepsilon(\sigma)} \Sym{w_0\lambda}.
\end{equation*}
Multiplying by $\Symo$ from the right then yields
\begin{equation*}
M_\lambda = \frac{q_{w_\lambda}^{-1}}{\W{}} \sum_{\mu \in X^\vee} q_{-\mu} \sum_{\substack{\sigma \in \Gamma_{t^\lambda}^+(\mu) \\ \varepsilon(\sigma) \in W^{w_0\lambda}}} q_{w_0\iota(\sigma)} L_\sigma X_{\mu} \Symo
\end{equation*}
and thus $L_{\lambda\mu} = q_{w_\lambda}^{-1} \sum_{\substack{\sigma \in \Gamma_{t^\lambda}^+(\mu) \\ \varepsilon(\sigma) \in W^{w_0\lambda}}} q_{w_0\iota(\sigma)} L_\sigma$. Moreover, we see that for a LS-gallery $\sigma$ we have $\varepsilon(\sigma) \in W^{w_0\lambda}$.
\end{remark}

\section{Commutation formula}
\label{sec:commutation}

Here we give two other applications of theorem~\ref{th:multiplicationgalleries}. The first is a nice combinatorial description of the base change matrix between the standard basis $\{T_w\}_{w \in \EWaffin}$ and the basis $\{X_\lambda \overline{T}_w\}_{\mu \in X^\vee, w \in W}$.

\begin{corollary}\label{co:standardasalcove}
Let $v \in \EWaffin$ and fix some minimal gallery of type $t$ connecting $A_f$ and $A_v$. Then
\begin{equation*}
T_v = \sum_{\sigma \in \Gamma^+_{t}, \iota(\sigma)= id} L_\sigma X_{e(\sigma)} 
= \sum_{\sigma \in \Gamma^+_{t}, \iota(\sigma) = id} q^{-1}_{wt(\sigma)} q_{\varepsilon(\sigma)} L_\sigma X_{wt(\sigma)} \overline{T}_{\varepsilon(\sigma)}.
\end{equation*}
\end{corollary}

Now let $w \in W^\lambda$ and recall the notation from definition~\ref{def:stabilizers}. Then $l(w w_\lambda) = l(w) + l(w_\lambda)$ and thus $T_{w w_\lambda} = T_w T_{w_\lambda}$. We also have $T_{\tau_\lambda} = T_{n^\lambda} T_{w_0} \overline{T}_{w_\lambda}$ by lemma~\ref{le:lengthfunction}. Since $T_{w_\lambda} X_\lambda = X_\lambda T_{w_\lambda}$ we get
\begin{equation*}
T_w T_{w_\lambda} X_\lambda = T_w X_\lambda T_{w_\lambda} = q_\lambda^{-1} T_w T_{n^\lambda} T_{w_0} \overline{T}_{w_\lambda} T_{w_\lambda} = q_\lambda^{-1} T_{w n^\lambda} T_{w_0}.
\end{equation*}
Applying corollary~\ref{co:standardasalcove} to $v = wn^\lambda$ and using $\overline{T}_y T_{w_0} = T_{y w_0}$ for $y \in W$ we get a $q$-analog of a commutation formula of Pittie and Ram in the nil-affine Hecke algebra~\cite{pittieram:99}.
\begin{corollary} 
Let $\lambda \in X^\vee_+$ and $w \in W^\lambda$. Let $t$ be the type of a minimal gallery connecting $A_f$ and $A_{wn^\lambda}$. Then
\begin{equation*}
T_{ww_\lambda} X_\lambda = q^{-1}_\lambda T_{wn^\lambda} T_{w_0} = \sum_{\sigma} q^{-1}_{wt(\sigma) +\lambda} q_{\varepsilon(\sigma)} L_\sigma X_{wt(\sigma)} T_{\varepsilon(\sigma)w_0}
\end{equation*}
where the sum is over all $\sigma \in \Gamma^+_t$ starting in $A_f$.
\end{corollary}

\section{Geometric interpretation}
\label{sec:geometricinterpretation}

In this section we reformulate results of~\cite{billigdyer:94} using galleries to show some relations of our combinatorics to geometry. We show that this geometrical interpretation is compatible with the geometrical interpretation of the $L_{\lambda\mu}$ arising from the geometric definition of the Satake isomorphism (see~\cite{haineskottwitzprassad:03}). Since there are no new results in this section we are rather sketchy.

Assume that $X^\vee = Q^\vee$. So $\EWaffin = \Waffin$ and generalized alcoves are alcoves. We use negatively folded galleries as in remark~\ref{re:negative}. We regard the affine Hecke algebra $\EHecke$ specialized at some prime power $\pq$ and all polynomials evaluated at $\pq$.

Details of the following constructions and their relation to affine Kac-Moody algebras can be found in Kumar's book~\cite{kumar:02}. We use the notation from the introduction, i.e. $G$ is a semisimple, simply connected algebraic group over $K$ with root datum $\Phi$ associated to a Borel $B \subset G$ and a maximal torus $T \subset B$ and all groups are defined and split over $\IF_\pq$. Let $\kk = K((t))$ be the field of Laurent series and denote by $\ko = K[[t]] \subset \kk$ the ring of formal power series. The evaluation at $0$ induces a map $ev : \ko \to K, t \mapsto 0$. This induces a morphism of groups $ev : G(\ko) \to G$. Define $\Kmg{B} = ev^{-1}(B)$. Further we set $\Kmg{G} = G(\kk)$ and let $\Kmg{N} \subset \Kmg{G}$ be the normalizer of $T$ in $\Kmg{G}$. Then $(\Kmg{G}, \Kmg{B}, \Kmg{N}, T(\ko))$ is a Tits system with Weyl group $\Waffin$.

From the theory of Tits system one has a Bruhat decomposition $\Kmg{G} = \bigsqcup_{w \in \Waffin} \Kmg{U}_w w \Kmg{B}$ with $\Kmg{U}_w$ isomorphic to $\IA_K^{l(w)}$. On the other hand there is the Iwasawa decomposition $\Kmg{G} = \bigsqcup_{w \in \Waffin} U(\kk) w \Kmg{B}$ where $U \subset B$ is the unipotent radical. These decompositions are compared by Billig and Dyer  in~\cite{billigdyer:94}.

\begin{theorem}[\cite{billigdyer:94}]\label{th:generalizedbruhat}
For $w \in \Waffin$ and $s = s_i \in \Saffin$ one has
\begin{equation*}
U(\kk) w \Kmg{B} s \Kmg{B} =
\begin{cases}
U(\kk) ws \Kmg{B} & \text{if } A_w \succ A_{ws}\\
U(\kk) w \Kmg{B} \sqcup U(\kk) ws \Kmg{B} & \text{if } A_w \prec A_{ws}.
\end{cases}
\end{equation*}
\end{theorem}

Let $w \in \Waffin$ and let $\sigma$ be a minimal gallery of type $t = (t_1, \ldots, t_k)$ which connects $A_f$ and $A_w$. Define a map $\eta :\Kmg{U}_w \to \Gamma^-_t$ by
\begin{equation*}
\eta(u) = (A_f, A_{w_1}, \ldots, A_{w_k}) \text{ iff  } u t_1 \cdot \ldots \cdot t_j \in U(\kk) w_j \Kmg{B} \text{ for } j \in \{1, \ldots, k\}.
\end{equation*}
It follows from the last theorem that $\eta$ is well defined.

The affine flag variety is the quotient $\Flaff$. It is the set of closed points of an ind-variety defined over $\IF_\pq$ such that the Bruhat cells $\Kmg{B} w \cdot \Kmg{B}$ for $w \in \Waffin$ are isomorphic to affine spaces of dimension $l(w)$. For $Y \subset \Kmg{G}$ denote by $Y \cdot \Kmg{B}$ its image in $\Flaff$. For a negatively folded gallery $\sigma$ denote by $m^-(\sigma)$ the total number of negative directions and by $n^-(\sigma)$ the total number of negative folds. Then one has

\begin{theorem}[\cite{billigdyer:94}]\label{th:galleriesbruhat}
\begin{enumerate}
\item If $\sigma \in \Gamma^-_t$ then $\eta^{-1}(\sigma) \cong K^{m^-(\sigma)} \times (K^*)^{n^-(\sigma)}$.
\item If $v \in \Waffin$ then $\Kmg{B} w \cdot \Kmg{B} \cap U(\kk) v \cdot \Kmg{B} = \bigsqcup_{\sigma \in \Gamma^-_{t} (A_f,A_v)} \eta^{-1}(\sigma) w \cdot \Kmg{B}$.
\end{enumerate}
\end{theorem}

For $Y \subset \Flaff$ denote by $|Y|$ the number of $\IF_\pq$-rational points. Let $L^-_w (A_f,A_v)$ be the analog of definition~\ref{de:lpolynomials} for negatively folded galleries. From the last theorem we get

\begin{corollary}\label{co:geometricinterpretation}
If $v,w \in \Waffin$ then $|\Kmg{B} w \cdot \Kmg{B} \cap U(\kk) v \cdot \Kmg{B}| = L^-_w (A_f,A_v)$.
\end{corollary}

\begin{remark}~\label{re:positiveopposite}
Looking at positively folded galleries starting in $-A_f$ one can calculate intersections $\Kmg{B}^- w \cdot \Kmg{B} \cap U^-(\kk) v \cdot \Kmg{B}$ in the same way as in corollary~\ref{co:geometricinterpretation}. Here $\Kmg{B}^-$ is obtained from the opposite Borel $B^-$ in the same way as $\Kmg{B}$ from $B$ and $U^-$ is the unipotent radical of $B$.
\end{remark}

As already mentioned in the introduction the definition of the \emph{geometric} Satake isomorphism yields a geometric interpretation of the $L_{\lambda\mu}$. For stating this we have to consider intersections in the affine Grassmanian $\Kmg{G}/G(\ko)$. For $\lambda, \mu \in X^\vee$ let $X_{\lambda\mu}^- \Def  G(\ko) \tau_\lambda \cdot G(\ko) \cap U^-(\kk) \tau_\mu \cdot G(\ko)$. Then $L_{\lambda\mu} = |X_{\lambda\mu}^-|$. This interpretation can also be recovered using the last corollary as follows.

The group $G(\ko)$ is the parabolic subgroup of  $\Kmg{G}$ associated to the classical Weyl group $W \subset \Waffin$, i.e. $G(\ko) = \bigsqcup_{w \in W} \Kmg{B} w \Kmg{B}$. Let $\pi : \Flaff \to \Kmg{G}/G(\ko)$ be the canonical projection. From general theory of Tits systems one knows that
\begin{equation*}
\pi_{|\Kmg{B} w \cdot \Kmg{B}} : \Kmg{B} w \cdot \Kmg{B} \to \Kmg{B} w \cdot G(\ko)
\end{equation*}
is an isomorphism iff $w$ is of minimal length in $wW$. So we get an isomorphism
\begin{equation*}
\pi_{|\sqcup_{w \in W^\lambda} \Kmg{B} w n^\lambda \cdot \Kmg{B}} :
\bigsqcup_{w \in W^\lambda} \Kmg{B} w n^\lambda \cdot \Kmg{B} \to
G(\ko) \tau_\lambda \cdot G(\ko).
\end{equation*}
If $\mu \in Q^\vee$ then
\begin{equation*}
\pi^{-1}(U(\kk) \tau_\mu \cdot G(\ko)) = \bigsqcup_{w \in W} U(\kk) \tau_\mu w \cdot \Kmg{B}.
\end{equation*}
Bringing this together and setting $X_{\lambda\mu} = G(\ko) \tau_\lambda \cdot G(\ko) \cap U(\kk) \tau_\mu \cdot G(\ko)$ we get
\begin{align*}
|X_{\lambda\mu}|
& =  \sum_{\substack{w \in W^\lambda \\ v \in W}} | \Kmg{B} w n^\lambda \cdot \Kmg{B} \cap U(\kk) \tau_\mu v \cdot \Kmg{B}|\\
& = \sum_{\substack{w \in W^\lambda \\ v \in W}} L^-_{w n^\lambda} (A_f,\tau_\mu v)
 = \sum_{w \in W^\lambda} L^-_{w n^\lambda} (A_f,\mu)\\
& = \frac{1}{W_\lambda(\pq)} L^-_{t^\lambda}(\mu) = \pq^{-2\langle \rho, \mu \rangle} L_{\lambda\mu}.
\end{align*}
The last equalities follow from the analog of the second statement in lemma~\ref{le:macdonaldasmonomial} for negatively folded galleries. Using the fact that $X^-_{\lambda\mu} = w_0 X_{\lambda, w_0\mu}$ and $L_{\lambda,\mu} = \pq^{\langle \rho, \mu-w_0\mu \rangle} L_{\lambda,w_0\mu} = \pq^{2\langle \rho, \mu \rangle} L_{\lambda,w_0\mu}$ we get as in~\cite{gaussentlittelmann:03} that
\begin{equation*}
|X^-_{\lambda\mu}| = \sum_{\sigma \in \Gamma^+_{t^\lambda}(\mu), \iota(\sigma) \in W^\lambda} \pq^{l(w_0\iota(\sigma))} L_\sigma.
\end{equation*}

It is well known that the affine Hecke algebra $\EHecke$ specialized at $\pq$ can be interpreted as the convolution algebra of $\Kmg{B}_\pq$-invariant functions with finite support on the affine flag  variety $\Kmg{G}_\pq / \Kmg{B}_\pq$ (see for example~\cite{haineskottwitzprassad:03}). In this setting the generator $T_w$ corresponds to the characteristic function on  $\Kmg{B}_\pq w \cdot \Kmg{B}_\pq$. Using theorem~\ref{th:generalizedbruhat} one can give a similar interpretation for the alcove basis.

Let $t_w \in F$ be the characteristic function on $U^-(\kk_\pq) w \cdot \Kmg{B}_\pq$ (which in general does not have finite support) and let $M$ be the subspace spanned by all $t_w$. Using theorem~\ref{th:generalizedbruhat} one can show that $M$ is closed under right convolution by $\EHecke$. Moreover, for $w \in \Waffin$ and $s \in \Saffin$ one gets
\begin{equation*}
t_w * T_s =
\begin{cases} 
t_{ws} & \text{if } A_w \prec A_{ws}\\
\pq t_{ws} + (\pq-1) t_w & \text{if } A_w \succ A_{ws}.
\end{cases}
\end{equation*}
So by lemma~\ref{le:qmultiplicationalcoves} the map
\begin{align*}
M & \to \EHecke\\
t_v & \mapsto \pq^{\langle \rho, \mu(A_v) \rangle} X_{wt(A_v)} \overline{T}_{\w(A_v)} = \pq^{2 \langle \rho, \mu(A_v) \rangle - l(\delta(A_v))} X_{A_v}
\end{align*}
is an isomorphism of right $\EHecke$-modules.

\begin{remark}
In~\cite{billigdyer:94} the cited results are shown for any Kac-Moody group and any generalized system of positive roots. Theorem~\ref{th:generalizedbruhat} is then formulated with distinguished subexpressions instead of folded galleries. All facts also follow quite immediately using the methods from~\cite{gaussentlittelmann:03} in this case.\\
The explicit formulas for the action of $\EHecke$ on $M$ are known (at least to experts). It is essentially the action of $\EHecke$ on the space of Iwahori fixed vectors of the universal unramified principal series. Corollary~\ref{co:geometricinterpretation} follows from these explicit formulas.
\end{remark}

\bibliography{bibliography}
\bibliographystyle{alpha}
\end{document}